\documentclass[11pt,a4paper]{amsart}
\usepackage{amssymb,amsmath,amscd,amsthm,latexsym,exscale}
\usepackage[all]{xy}

\def\cal{\mathcal}
\newcommand{\maxz}{\mbox{Max$\mathfrak Z$ }}
\newcommand{\f}[1]{\mbox{$\mathfrak #1$}}
\newcommand{\gmod}{\mbox{$\f{g}$-mod}}
\newcommand{\verma}[1]{\mbox{$M(#1)$}}
\newcommand{\ein}[2]{\mbox{${#1}_{\rvert_{\scriptscriptstyle {#2}}}$}}

\newcommand{\Unull}{\cal U_0}
\newcommand{\Uhmu}{\cal U_{\{h_{\mu}\}}}
\newcommand{\U}{\cal U}

\newcommand{\N}{{\cal N}_{\nu}}
\newcommand{\Ha}{{\cal H}_{\nu}}
\newcommand{\Seins}{{\cal S}_1}
\newcommand{\Szwei}{{\cal S}_2}

\newcommand{\hotimes}{\hat{\otimes}}
\newcommand{\F}{\phi}

\DeclareMathOperator{\pr}{pr}
\DeclareMathOperator{\fpr}{\bar{pr}}
\DeclareMathOperator{\id}{id}
\DeclareMathOperator{\Id}{Id}
\DeclareMathOperator{\sym}{sym}
\DeclareMathOperator{\PO}{PO}
\DeclareMathOperator{\p}{pro}
\DeclareMathOperator{\Quot}{Quot}
\DeclareMathOperator{\Ext}{Ext}
\DeclareMathOperator{\End}{End}
\DeclareMathOperator{\Hom}{Hom}
\DeclareMathOperator{\Homtau}{Hom_{{\cal M}(\chi(\tau))\to}}
\DeclareMathOperator{\Homtauinf}{Hom_{{\cal M}^{\infty}(\chi(\tau))\to}}
\DeclareMathOperator{\Ann}{Ann}
\DeclareMathOperator{\dif}{d}
\DeclareMathOperator{\codim}{codim}
\DeclareMathOperator{\can}{can}
\DeclareMathOperator{\nat}{nat}
\DeclareMathOperator{\Nat}{Nat}
\DeclareMathOperator{\tr}{tr}
\DeclareMathOperator{\inc}{i}
\DeclareMathOperator{\cont}{cont}
\DeclareMathOperator{\jott}{j}
\DeclareMathOperator{\ce}{c}
\DeclareMathOperator{\adj}{adj}

\DeclareMathOperator{\bDelta}{\Diamond}
\DeclareMathOperator{\Deltag}{\Delta_g}

\newtheorem{theo}{Theorem}
\newtheorem*{theoohne}{Theorem}
\newtheorem*{proposition}{Proposition} 
\newtheorem{lemma}{Lemma}
\newtheorem*{kor}{Corollary}
\newtheorem*{korbew}{Corollary of the proof}
\newtheorem*{norm}{Normalization identities}
\newtheorem*{zer}{Decomposition identity}
\newtheorem*{rot}{Rotation identity}
\newtheorem*{flach}{Flat triangle identity}
\newtheorem*{Split}{Split identities}
\newtheorem*{bei}{{\bf Example}}
\theoremstyle{remark}
\newtheorem*{bemerkung}{Remark}

\begin{document}
     

\title{The fine structure of translation functors}
\author{Karen G\"unzl}
\address{Mathematisches Institut\\ Universit\"at Freiburg\\ 
                           Eckerstr.1 \\ 79104 Freiburg \\ FRG}
\email{karen@mathematik.uni-freiburg.de}
\date{August 1998}
\begin{abstract}
   Let $E$ be a simple finite dimensional representation of a
      semisimple Lie algebra with extremal weight $\nu$ and let 
      $0 \neq e \in E_{\nu}$.
   Let $\verma{\tau}$ be the Verma module with highest weight $\tau$ and
       $0 \neq v_{\tau} \in \verma{\tau}_{\tau}$.
   We investigate the projection of $e \otimes v_{\tau} \in E
       \otimes \verma{\tau}$ on the central character $\chi(\tau+\nu)$.
   This is a rational function in $\tau$ and we calculate its poles and
       zeros.
   We then apply this result in order to compare translation functors.
\end{abstract}
\maketitle

     \tableofcontents 
     \section{\bf Introduction}

Let $k$ be a field of characteristic zero and let $\f{g}$ be a split
   semisimple Lie algebra over $k$.
Then 
   $\f{h} \subset \f{b} \subset \f{g}$ denote a Cartan and a Borel subalgebra, 
   $\f{U}= \f{U}(\f{g})$ the universal enveloping algebra and
   $\f{Z} \subset \f{U}$ its center. 
We set $R^+ \subset R \subset \f{h}^*$ to be the roots of \f{b} and of \f{g}.
Let \f{n} (resp. \f{n^-}) be the sum of the positive (negative) weight 
   spaces and denote by $P^+ \subset P \subset \f{h}^*$ the set of 
   dominant and the lattice of integral weights respectively.
Let $\cal W$ be the Weyl group.

A $\f{g}$-module $M$ is called {\em$\f{Z}$-finite}, 
                   if $\dim_k (\f{Z} m) < \infty$ for all $m \in M$.
Every $\f{Z}$-finite module $M$ splits under the operation of the center 
   $\f{Z}$ into a direct sum of submodules 
   $M = \oplus_{\chi \in \maxz}M_{\chi}$.
Here  $\chi$ runs over the maximal ideals in the center and 
   $M_{\chi} \in {\cal M}^{\infty}(\chi)$, where
   \[ {\cal M}^{\infty} (\chi) := \{ M  \;\mid   
      \mbox{ for all } m \in M \mbox{ there exists } n \in \mathbb N 
      \mbox{ such that } \chi^n m = 0 \}.\]

We denote the projection onto the central character $\chi$ by 
   $\pr_{\chi}:\, M \mapsto M_{\chi}$. 
By means of the  Harish-Chandra homomorphism \cite[7.4]{di} 
   $\xi:\,\f{Z} \to S(\f{h})$ (normalized by
   $z-\xi(z) \in \f{U}\f{n}\; \forall z \in \f{Z}$)
   we assign to each weight $\tau \in \f{h}^*$ its central character
   $\chi(\tau)$ defined by $\chi_{\tau}(z) := \big(\xi(z)\big)(\tau) \;
   \forall z \in \f{Z}$.

Let now $\nu \in P$ be an integral weight and denote by $E :=E(\nu)$ 
   the irreducible finite dimensional $\f{g}$-module with highest
   weight in $\cal W \nu$.
It is known that for each $\f{Z}$-finite module $M$ the tensor product
   $E \otimes M$ is again $\f{Z}$-finite \cite[Corollary 2.6]{bg}.
If we now tensor a module $M \in {\cal M}^{\infty} \big(\chi(\tau)\big)$
with $E$ and then project onto the central character
$\chi(\tau +\nu)$, we obtain the {\em translation functor}
  {\renewcommand{\arraystretch}{1.5}
  \[
  \begin{array}{lccc}
     T_{\tau}^{\tau +\nu}:\, & \cal M^{\infty} \big(\chi(\tau)\big)& 
     \longrightarrow & \cal M^{\infty} \big(\chi(\tau+\nu)\big)\\
     &   M   & \mapsto & \pr_{\chi(\tau+\nu)}(E \otimes M)
  \end{array}
  \]}

In the present paper we investigate the fine structure of translation functors.
Namely, take for $M$ the Verma module
$\verma{\tau}:= \f{U}(\f{g}) \otimes_{\f{U}(\f{b})} k_{\tau}$ with its highest
weight vector
$v_{\tau} := 1 \otimes 1 \; \in \verma{\tau}_{\tau}$ and then identify for all
$\tau \in \f{h}^*$:
   \[
   \begin{array}{ccc}
       \verma{\tau} & \stackrel{\sim}{\longrightarrow} & \f{U}(\f{n}^-)\\
          u v_{\tau}&  \mapsto                         & u
   \end{array}
   \]

We set $V(\nu)= V:= E \otimes \f{U}(\f{n}^-)$, choose a fixed extremal 
   weight vector $e_{\nu} \in E_{\nu}$ and define the map
             \[
             \begin{array}{lccc}
              f_\nu:\, &\f{h}^* & \longrightarrow & V\\
              & \tau    & \mapsto         & \pr_{\chi(\tau+\nu)}
                                                  (e_{\nu} \otimes v_{\tau})
             \end{array}
             \] 
Here we identify $\pr_{\chi(\tau+\nu)}(e_{\nu} \otimes v_{\tau}) \in
   E \otimes \verma{\tau}$ with its image in $V$. 
The image of the map $f_{\nu}$ is then contained in the finite dimensional
   $\nu$-weight space $(E \otimes \f{U}(\f{n}^-))_{\nu}$ of $V$, 
   and we may thus regard $f_{\nu}$ as a map between varieties.
In general however, $f_{\nu}$ is not a morphism.

Let $\rho := 1/2 \sum_{\alpha \in R^+} \alpha$ denote the half sum of 
   positive roots, $\alpha^{\vee}$ the co-root of $\alpha$ and set
   \[\N :=\{(\alpha, m_{\alpha}) \in R^+ \times \mathbb Z \mid
          -\langle \rho, \alpha^{\vee} \rangle \le m_{\alpha}
             < - \langle \nu +\rho, \alpha^{\vee} \rangle \}.\]
Then the set
   \[\Ha := \bigcup_{(\alpha, m) \in \N}
         \{ \tau \in \f{h}^* \mid \langle \tau,\alpha^{\vee} \rangle = m \} \] 
   is a finite family of hyperplanes and therefore Zariski closed. 
We will show at first that $f_{\nu}$ is a morphism of varieties 
   on the complement of $\Ha \cup \cal S$, where $\cal S \subset \f{h}^*$ 
   is a suitable Zariski closed subset of codimension $\ge 2$.
More precisely, $\cal S$ consists of intersections of finitely many
   hyperplanes.

Define for all roots $\alpha \in R$ and for all $m \in \mathbb Z$
   the polynomial function $H_{\alpha,m}$ on $\f{h}^*$ by 
   \[H_{\alpha, m}(\tau) := \langle \tau,\alpha^{\vee} \rangle -m 
                         \mbox{ for all } \tau \in \f{h}^*.\]
If we then set
   $\delta_{\nu} :=\prod_{(\alpha, m) \in \N} H_{\alpha, m}$,
   we obtain $\Ha = \{\tau \in \f{h}^* | \; \delta_{\nu}(\tau) =0 \}$.
A central result of this paper will be the following
\begin{theoohne}
   There exists a morphism of varieties 
   $G : \, \f{h}^* \to V_{\nu}$, such that the set of zeros of $G$ has 
   codimension $\ge 2$ and such that $G$ equals $\delta_{\nu} f_{\nu}$
   on $\f{h}^* - (\Ha \cup \cal S)$.
\end{theoohne}

This means that the map $f_{\nu}$ has a pole of order $1$ along each 
   of the hyperplanes 
   $\langle \tau, \alpha^{\vee} \rangle = m$ with $(\alpha, m) \in \N$
   and outside of these hyperplanes it is a non-vanishing morphism of 
   varieties except on a set of codimension $\ge 2$.\\

In chapter \ref{kapvier} we introduce the so-called triangle functions
   $\Delta(\mu,\nu;x)(\tau)$ for integral weights $\mu$ and $\nu$ and
   $x \in \cal W$. 
These are rational functions on $\f{h}^*$, which measure in a subtle way 
   the relation between the two translation functors
   $T_{\tau+\nu}^{\tau+\nu+\mu} \circ T_{\tau}^{\tau+\nu}$ and 
   $T_{\tau}^{\tau+\nu+\mu}$ by first applying them to the Verma modules
   $\verma{\tau}$ and $\verma{x\cdot\tau}$ and then identifying the results
   with $\verma{\tau+\nu+\mu}$ and $\verma{x\cdot(\tau+\nu+\mu)}$
   respectively.
   
 Bernstein defined in \cite{be} the so-called {\em relative trace $\tr_E$} 
     for a finite dimensional vector space $E$ and this function is related
     to the special case $\Delta(-\nu,\nu;w_0)$, where $\nu$ is dominant and 
     $w_0$ is the longest element of $\cal W$.
 This is explained in Chapter \ref{bern}.
 By Bernstein's explicit formula for $\tr_E$ we obtain

   \begin{kor}
      Let $\nu \in P^+$ be dominant and $\tau \in \f{h}^*$ generic, i.e.
      $\langle \tau, \alpha^{\vee}\rangle \notin \mathbb Z$ for all
      $\alpha \in R$.
      Then
      \[\Delta(-\nu, \nu;w_0)(\tau) 
      =
      \prod_{\alpha \in R^+} \frac{\langle \tau +\nu +\rho, 
                                             \alpha^{\vee} \rangle}
                                  {\langle \tau +\rho, 
                                             \alpha^{\vee} \rangle} \; .\]
   \end{kor}
 
 We remark that Kashiwara has also calculated this case \cite[Thm. 1.9]{ka}.

In Chapter \ref{formel} we calculate the triangle functions in general.
Since they are rational functions on $\f{h}^*$, it suffices to determine their
    zeros and poles.
In order to do this we make use of the maps $f_{\nu}$ for suitable integral
    weights $\nu$.
Set now $\bar{\alpha}(\lambda) := 1$ 
                for $\langle\lambda,\alpha^{\vee}\rangle < 0$ and
   $\bar{\alpha}(\lambda) := 0$ 
                for $\langle\lambda,\alpha^{\vee}\rangle \ge 0$.
Then we obtain

   \begin{theoohne}
      Let $\nu, \mu \in P$ be integral weights and $x \in \cal W$. 
      Then
     \[\Delta(\mu, \nu; x)(\tau-\rho)=
       c \!\!\!\!\!\!\prod_{\genfrac{}{}{0pt}{2}
                {\alpha \in R^+}{\text{with } x\alpha \in R^-}}\!\!\!\!
       \frac
        {\langle \tau,\alpha^{\vee}\rangle^{\bar{\alpha}(\nu)}
         \langle\tau+\nu,\alpha^{\vee}\rangle^{\bar{\alpha}(\mu)}
         \langle\tau+\nu+\mu,\alpha^{\vee}\rangle^{\bar{\alpha}(\nu+\mu)}}
        {\langle \tau,\alpha^{\vee}\rangle^{\bar{\alpha}(\nu+\mu)}
         \langle \tau+\nu,\alpha^{\vee}\rangle^{\bar{\alpha}(\nu)}
         \langle \tau+\nu+\mu,\alpha^{\vee}\rangle^{\bar{\alpha}(\mu)}}\]
        for a constant $c \in k^{\times}$ independent of $\tau, \nu, \mu$ 
        and $x$.
   \end{theoohne} 

\subsubsection*{\bf Acknowledgment}
   I am grateful to my advisor, Wolfgang Soergel, for his constant 
   encouragement, his patience and all the fruitful discussions we had.
   Without his many ideas this work would not have been possible.
   
   The author was partially supported by EEC TMR-Network ERB FMRX-CT97-0100.

     \section{\bf Filtration of $\boldsymbol{E \otimes \verma{\tau}}$}

In the following let $\nu \in P$ be a fixed integral weight.
Then let $\nu_1,\ldots,\nu_n$ be the multiset of weights of $E=E(\nu)$,
     i.e. with multiplicities, 
     and let $(e_i)_{1\le i \le n}$ be a basis of weight vectors of $E$,
     such that $e_i \in E_{\nu_i}$ and $\nu_i < \nu_j \Rightarrow i < j$.
Here $\lambda \le \mu$ for weights $\lambda, \mu \in \f{h}^*$, if 
     $\mu-\lambda = \sum_{\alpha \in R^+}n_{\alpha}\alpha$ with 
     $n_{\alpha} \in \mathbb N$.

In \cite{bgg1} it was shown that the tensor product $E \otimes \verma{\tau}$
     admits a chain of submodules 
     $N_i := \sum^n_{j=i} \f{U}(\f{n}^-) (e_j \otimes v_{\tau})$, 
     such that $N_i / N_{i+1} \cong \verma{\tau + \nu_i}$.
 
Using this, one can easily construct a slightly coarser filtration, namely a
     filtration such that the subquotients are each isomorphic to a direct sum
     of Verma modules with the same highest weight.
In order to do this, let $\mu_1,\mu_2,\ldots,\mu_m$ denote the set of weights
     of $E$ {\em without} multiplicities and set $d_j := \dim E_{\mu_j}$.
We may choose the numbering of the $\nu_i$ in such a way that 
     $\nu_k = \mu_j$ for all $k$ with 
     $d_1 + \cdots +d_{j-1} < k \le d_1 + \cdots +d_j$.
The weight space $E_{\mu_j}$ is then generated by precisely these $e_k$.
Set $d_0 :=0$ and define for $1\le j \le m$:
      \[M_j :=\sum\limits^m_{k=j} \f{U}(\f{n}^-) (E_{\mu_k} \otimes v_{\tau})
             =\sum\limits^n_{k>d_1+d_2+\cdots d_{j-1}} 
                              \f{U}(\f{n}^-) (e_k \otimes v_{\tau}) .\]

We obtain immediately: \label{filtration}
The chain $E \otimes \verma{\tau} =M_1 \supset \cdots \supset M_m \supset 0$
     is a filtration of submodules such that 
     $M_j /M_{j+1} \cong {\bigoplus}_{d_j} \verma{\tau +\mu_j}$
and the summands $\verma{\tau +\mu_j}$ are generated by the vectors
     $(e_i \otimes v_{\tau})/M_{j+1}=
     \pr(e_i \otimes v_{\tau})/M_{j+1}$, where $e_i$ is a vector of weight
     $\mu_j$.

\begin{lemma} \label{lemma_Z}
       For $1 \le j \le m$ let $s_j \in \f{Z}$ with
       $\chi_{\tau +\mu_j}(s_j) =0$ and set $z_i :=\prod_{j \ge i} s_j$.
       Then for $1 \le i \le m$ we get: $z_i M_i = 0$.

       In particular, for a weight vector $e_{\mu_i} \in E_{\mu_i}$ we obtain
       \[z_{i+1} (e_{\mu_i} \otimes v_{\tau})=
            z_{i+1} \pr_{\chi(\tau +\mu_i)}(e_{\mu_i} \otimes v_{\tau}).\]
\end{lemma}

\begin{proof}
     The first statement follows by induction from above and since a central
     element $z \in \f{Z}$ operates on the Verma module $\verma{\lambda}$ by
     multiplication with the scalar $\chi_{\lambda}(z)$.
 
     By construction of the $M_i$ it is clear, that for
     $e_{\mu_i}\otimes v_{\tau} \in M_i$ its projections 
     $\pr_{\chi(\tau +\mu)}(e_{\mu_i} \otimes v_{\tau})$ lie already in
     $M_{i+1}$ if $\chi(\tau +\mu) \neq \chi(\tau +\mu_i)$.
     But $e_{\mu_i}\otimes v_{\tau}$ equals the sum of its projections and
     thus the second statement follows from the first.
\end{proof}

Note that for the extremal weight $\nu = \nu_{i_0} = \mu_{j_0}$ 
     its weight space is of dimension $1$.
Therefore, we can always choose the numbering of the $\nu_i$ 
     such that $i>i_0 \iff \nu_i > \nu_{i_0}=\nu$ or equivalently 
     $j>j_0 \iff \mu_j > \mu_{j_0}=\nu$.
It then follows immediately that an element of $M_{i_0} = M_{j_0}$ 
     (of $M_{i_0+1} = M_{j_0+1}$ resp.) consists only of parts with central 
     character $\chi_{\tau +\mu}$ for $\mu \in P(E)$ with $\mu \ge \nu$
     ($\mu > \nu$ resp.).  
For $e_{\nu} \otimes v_{\tau} \in M_{i_0}$ we thus obtain the following
     special case of Lemma \ref{lemma_Z}

\begin{lemma} \label{spezial}
       For $\mu \in P(E)$ with $\mu> \nu$ let $s_{\mu} \in \f{Z}$
       such that $\chi_{\tau +\mu}(s_{\mu}) =0$ and
       set $z:= \prod_{\mu >\nu} s_{\mu}$. Then
       $z (e_{\nu} \otimes v_{\tau})=
            z \pr_{\chi(\tau +\nu)}(e_{\nu} \otimes v_{\tau}).$  
\end{lemma}

     \pagebreak
\section{\bf The fine structure of $\boldsymbol{E \otimes M(\tau)}$}

\subsection{Algebraicity of $f_{\nu}$} 

  For $\nu \in P$ fixed, set $E := E(\nu)$ and let $f_{\nu}$ be the map
     \[
     \begin{array}{lccc}
     f_\nu:\, &\f{h}^* & \longrightarrow & V\\
              & \tau    & \mapsto        & \pr_{\chi(\tau+\nu)}
                                                  (e_{\nu} \otimes v_{\tau})
     \end{array}
     \] 
    where $V= V(\nu):= E \otimes \f{U}(\f{n}^-) 
          \cong E \otimes \verma{\tau}$ for all $\tau \in \f{h}^*$.
We will first construct a Zariski open set $\U \subset \f{h}^*$ such that
   $f_{\nu}$ restricted to $\U$ is a morphism of varieties.
In this case we will also say that $f_{\nu}$ is {\em algebraic} on $\U$.

Before we go on, we need some more notation : 
The dot-operation of the Weyl group on $\f{h}^*$ has fixed point 
   $-\rho = -1/2 \sum_{\alpha \in R^+} \alpha$ and is defined by
   $w \cdot \lambda := w(\lambda+\rho)-\rho$.
Let $\cal P (\f{h}^*)$ denote the set of polynomial functions on $\f{h}^*$
   and $\cal P (\f{h}^*)^{\cal W \cdot}$ the $(\cal W \cdot)$-invariants.
We define the operator
   \[\sym :\,
                \cal P (\f{h}^*) \longrightarrow 
                \cal P (\f{h}^*)^{\cal W \cdot}\]
   by $(\sym s)(\lambda)
             := \prod\limits_{x \in \cal W} s(x \cdot \lambda)$ \;\;
             for all $\lambda \in \f{h}^* , s \in \cal P (\f{h}^*)$.

For any $\mu \in P(E)$ with $\mu >\nu$ we now choose an element
   $h_{\mu} \in \f{h}$, $h_{\mu} \neq 0$ and define maps
   $H_{\mu}: \, \f{h}^* \times \f{h}^* \longrightarrow k$ by
   \[H_{\mu}(\lambda,\tau) := \langle \lambda - \tau -\mu,h_{\mu} \rangle
                         \;\;\; \mbox{ for all } \lambda, \tau \in \f{h}^*.\]

If we fix $\tau \in \f{h}^*$, the map 
   $H_{\mu}(-,\tau) :\, \f{h}^* \longrightarrow k$
   is then a polynomial function on $\f{h}^*$ and its kernel is
   obviously a hyperplane in $\f{h}^*$. 
We then define a $({\cal W \cdot})$-invariant polynomial function $p_{\tau}$,
   depending also on the choice of the elements $h_{\mu}$, by
\[p_{\tau}:= \prod_{\mu \in P(E), \mu >\nu} \sym H_{\mu}(-,\tau).\]

By means of the Harish-Chandra isomorphism \cite[7.4]{di} 
     $\xi:\,\f{Z} \stackrel{\sim}{\longrightarrow}
                                 {\cal P}(\f{h}^*)^{\cal W \cdot}$,  
     there exist central elements $s_{\mu} \in \f{Z}$, such that
     $\xi(s_{\mu}) =\sym H_{\mu}(-,\tau)$.
For $z_{\tau} := \prod_{\mu >\nu} s_{\mu}$ we have
     $\xi(z_{\tau})= p_{\tau}$.
Define now the map 
             \[
             \begin{array}{lccc}
             u=u_{\{h_{\mu}\}}:\, & \f{h}^* & \longrightarrow & k \\
                  & \tau    & \mapsto       & p_{\tau}(\tau +\nu)  
             \end{array}
             \]
and set 
    $\Uhmu :=\{ \tau \in \f{h}^* \mid  u(\tau) \neq 0 \} \subset \f{h}^*$. 

\begin{samepage} 
\begin{lemma}  \label{resultat}
   The map $f_{\nu}$ restricted to \; $\Uhmu$ is a morphism of varieties and
   for all $\tau \in \Uhmu$ we have 
   $f_{\nu}(\tau) = 
      \left(u(\tau)\right)^{-1} z_{\tau} (e_{\nu} \otimes v_{\tau})$.
\end{lemma}
\end{samepage}
\begin{proof}
   We know $\chi_{\tau +\mu}(s_{\mu})
                     = (\xi(s_{\mu}))(\tau +\mu) 
                     = (\sym H_{\mu}(-,\tau))(\tau+\mu)
                     = 0$
   and by Lemma \ref{spezial} we conclude 
               $z_{\tau}(e_{\nu} \otimes v_{\tau}) 
                = z_{\tau} \pr_{\chi(\tau+\nu)}(e_{\nu} \otimes v_{\tau})$.
   
   Let now $\tau \in \Uhmu$, i. e. $u(\tau)\neq 0$. 
   This means, that for all
   $w \in \cal W$ and for all $\mu \in P(E)$ with $\mu >\nu$ we have
   $\langle w \cdot (\tau+\nu)-\tau-\mu, h_{\mu} \rangle \neq 0$ and in
   particular $w \cdot (\tau+\nu) \neq \tau+\mu$ or equivalently
   $\chi(\tau +\nu) \neq \chi(\tau +\mu) \forall \mu \in P(E), \mu >\nu$.
   This implies, that already the projection
   $\pr_{\chi(\tau +\nu)}(e_{\nu} \otimes v_{\tau})
   \in E \otimes \verma{\tau}$ generates a Verma module 
   $\verma{\tau +\nu}$. 
   A central element $z$ operates on this by multiplication 
   with the scalar $(\xi(z))(\tau+\nu) = u(\tau)$. 
   We therefore get
           $ z_{\tau}(e_{\nu} \otimes v_{\tau}) 
              = z_{\tau} \pr_{\chi(\tau+\nu)}(e_{\nu} \otimes v_{\tau}) 
              = u(\tau) f_{\nu}(\tau)$ 
   and the above equation follows.

   We have yet to show that $f_{\nu}$ is algebraic on $\Uhmu$.
   By construction of $z_{\tau}$ it is clear that
   $z_{\tau}$ depends algebraically on $\tau$ for all $\tau \in \f{h}^*$. 
   Also the map
   $\tau \mapsto z_{\tau} (e_{\nu} \otimes v_{\tau}) \in V$ is
   a morphism on $\f{h}^*$ and $1/u$ is per definition algebraic
   on $\Uhmu$. Thus $f_{\nu}$ restricted to $\Uhmu$ is a morphism.
\end{proof}


\begin{bei} \label{beispiel}
   Let $\nu \in P(E)$ be an extremal and dominant weight.
   For all $\mu \in P(E)$ with $\mu \neq \nu$ we then have $\mu < \nu$
   and we may choose $z_{\tau} := 1 \in \f{Z}$ for all $\tau \in \f{h}^*$. 
   The above lemma then implies $\Uhmu = \f{h}^*$, the map $f_{\nu}$
   is a morphism on $\f{h}^*$ and 
   $f_{\nu}(\tau)= e_{\nu} \otimes v_{\tau}$ for all $\tau$. 
\end{bei}

  Lemma \ref{resultat} implies that $f_{\nu}$ is algebraic on 
    ${\U} := \bigcup \Uhmu$, where the union is taken over all possible
    choices for $\{h_{\mu} \mid \mu \in P(E), \mu > \nu \}$.
  If we set   
    \[{\cal A} := \{ \tau \in \f{h}^* \mid u(\tau)= p_{\tau}(\tau +\nu) =0
                    \mbox{ for all choices of } \{h_{\mu}\}  \} ,\]
    we can write $\U$ as $\U = \f{h}^* - \cal A$.
  Since $p_{\tau}(\tau +\nu) = 
         \prod_{\mu \in P(E), \mu > \nu} \prod_{w \in \cal W}
         \langle w \cdot (\tau +\nu)- \tau - \mu, h_{\mu} \rangle$
    we know that $\tau$ is in $\cal A$ if and only if there exists a
    $w \in \cal W$ and a $\mu \in P(E)$ with $\mu > \nu$, such that
    $w \cdot (\tau +\nu)= \tau + \mu$.
  Thus we obtain
    \[{\cal A} =\bigcup_{\genfrac{}{}{0pt}{2}
                        {\mu \in P(E), \mu >\nu}{w \in \cal W}}
                \{ \tau \in \f{h}^* \mid w \cdot (\tau +\nu) = \tau + \mu\}. \]
  Let us examing more closely the sets on the right hand side:
  For $w = s_{\alpha}, \alpha \in R$, we get 
    $\{ \tau \in \f{h}^* \mid \langle \tau, \alpha^{\vee}\rangle \alpha =
          \nu -\mu-\langle \nu +\rho, \alpha^{\vee} \rangle \alpha \}$,
    which implies that this set is nonempty if and only if there exists
    a weight $\mu = \nu-\langle \tau +\nu +\rho, \alpha^{\vee} \rangle \alpha$
    in the $\alpha$-string through $\nu$, such that $\mu \in P(E)$ and 
    $\mu > \nu$.   
  Since $\nu$ was an extremal weight, it is either the greatest or the
    smallest weight in this $\alpha$-string and hence such a $\mu$ exists
    only if $\langle \nu, \alpha^{\vee} \rangle < 0$ for an $\alpha \in R^+$.
  In this case all $\mu_{(n)} := \nu + n \alpha$ for 
    $1 \le n \le -\langle \nu, \alpha^{\vee} \rangle$ are weights of $E$ 
    with $\mu_{(n)} \ge \nu$.     
  Comparing this with 
    $\mu = \nu - \langle \nu + \rho + \tau, \alpha^{\vee} \rangle \alpha$
    we obtain as condition for $\tau$ precisely
    $-\langle \rho,\alpha^{\vee} \rangle \le
    \langle \tau, \alpha^{\vee} \rangle \le 
    -\langle \nu+\rho, \alpha^{\vee} \rangle - 1$ for an $\alpha \in R^+$.
  
Set ${\cal W}^0:=\{w \in {\cal W} \mid w \neq s_{\alpha} {\text{ for an }} 
                                          \alpha \in R \}$.
For
    \begin{eqnarray*}
    \N  &:=& \{(\alpha, m_{\alpha}) \in R^+ \times \mathbb Z \mid
            -\langle \rho, \alpha^{\vee} \rangle \le m_{\alpha}
             < - \langle \nu +\rho, \alpha^{\vee} \rangle  \} \\
    \Ha &:=& \bigcup_{(\alpha, m) \in \N}
             \{\tau \in \f{h}^* \mid\langle\tau, \alpha^{\vee}\rangle = m \}\\ 
    {\cal S} &:=& 
          \bigcup_{\genfrac{}{}{0pt}{2}
                  {\mu \in P(E), \mu > \nu}{w \in {\cal W}^0}} 
                  \{\tau \in \f{h}^* \mid w \cdot (\tau +\nu) = \tau +\mu \} ,
    \end{eqnarray*}

    we get ${\cal A} = \Ha \cup \cal S$ and $\cal S$ consists of finitely many
    intersections of at least two hyperplanes and is therefore a Zariski 
    closed subset of codimension $\ge 2$.
  We remark that in general $\Ha$ and $\cal S$ are not disjoint.  
  Using Lemma \ref{resultat} we obtain
\begin{samepage}
\begin{lemma} \label{algebraisch}
   The map $f_{\nu}$ is algebraic on the complement of $\Ha \cup \cal S$. 
\end{lemma}
\end{samepage}


\subsection{Poles of $f_{\nu}$. Theorem \ref{theoeins} and proof}\label{pole} 

Let now $\delta_{\nu} \in {\cal P}(\f{h}^*)$ be the product of the 
    equations of the hyperplanes in $\Ha$, that is 
    \[ \delta_{\nu} := \prod_{(\alpha, m) \in \N} H_{\alpha, m} 
                     =\prod_{\alpha\in R^+}\;\;
                  \prod_{0\le n_{\alpha}< -\langle\nu, \alpha^{\vee}\rangle}
                  (\langle \tau +\rho,\alpha^{\vee} \rangle -n_{\alpha})\: ,
    \] 
   where $H_{\alpha, m}(\tau) := \langle\tau, \alpha^{\vee} \rangle - m$. 
   We then have $\Ha = \{\tau \in \f{h}^* | \; \delta_{\nu}(\tau) = 0 \}$. 
\begin{samepage}
\begin{theo} \label{theoeins}
   There exists a morphism of varieties 
   $G : \, \f{h}^*  \longrightarrow V_{\nu}$,
   such that the set of zeros of $G$ has codimension $\ge 2$
   and such that $G$ equals $\delta_{\nu} f_{\nu}$
   on $\f{h}^* - (\Ha \cup \cal S)$.
\end{theo}
\end{samepage}


   The proof comes in two parts. In the first part we demonstrate the 
   existence of an algebraic extension $G$ of $\delta_{\nu} f_{\nu}$ on
   the whole of $\f{h}^*$, in the second we will show that the set of
   zeros of $G$ has codimension $\ge 2$.
   It will be useful to introduce the set $\Seins \subset \f{h}^*$ of all
   intersections of hyperplanes in $\Ha$. So if we set
   \[\Seins :=\{\tau\in\f{h}^* \mid H_{\alpha, m}(\tau) = 0 = 
      H_{\beta, n}(\tau) \mbox{ for } (\alpha, m) \neq (\beta, n) \in \N \} \]
   then ${\cal S} \cup \Seins$ is a Zariski closed subset of codimension 
   $\ge 2$.
   Note that $(\Ha \cap {\cal S}) \subset \Seins$ and hence
   $\Ha - ({\cal S} \cup \Seins) = \Ha - \Seins$.
   By definition of $\Seins$ we obtain immediately
   \begin{lemma} \label{salpha}
        Let $\tau \in \Ha - \Seins$. Then there exists exactly one
        $\alpha \in R^+$, such that
        $s_{\alpha} \cdot (\tau +\nu) = \tau +\mu_0$ for a weight
        $\mu_0$ of $E$ with $\mu_0 >\nu$, namely 
        $\mu_0 = \nu - \langle \tau +\nu +\rho, \alpha^{\vee} \rangle \alpha$.
        For all $\mu \in P(E)$ with $\mu > \nu$ and $\mu \neq \mu_{0}$
        we have $w \cdot (\tau +\nu) \neq \tau +\mu$ for all $w \in \cal W$.
   \end{lemma}

 
  \begin{proof}[Proof of Theorem \ref{theoeins}]
        Again let $\U = \f{h}^* -(\Ha \cup \cal S )$ be the set 
        on which $f_{\nu}$ is algebraic. First we claim
  \begin{itemize}
    \item[($\ast$)] For all $\tau_0 \in \f{h}^* - ({\cal S} \cup \Seins)$ 
           there exists a Zariski open neighborhood $\Unull$
           of $\tau_0$ and an algebraic map
           $G_0 : \, \Unull \to V_{\nu}$,
           such that $G_0 =\delta_{\nu} f_{\nu}$ on $\U \cap \Unull$.
  \end{itemize}
 
   This statement then implies the existence of an algebraic extension
   of $\delta_{\nu} f_{\nu}$ {\em locally} around every point of 
   $\f{h}^* -({\cal S} \cup \Seins )$.
   This extension is unique and thus we obtain a {\em global} algebraic 
   extension on $\f{h}^* -({\cal S} \cup \Seins)$. 
   Since $\codim ({\cal S} \cup \Seins) \ge 2$ we can extend this to a
   morphism $G$ on the whole of $\f{h}^*$. 

   Let now $\tau_0 \in \f{h}^* - ({\cal S} \cup \Seins)$.

   If $\tau_0 \notin \Ha$, then 
   $\tau_0 \in \U = \f{h}^* - ( \Ha \cup \cal S )$ and since $f_{\nu}$ was
   algebraic on $\U$ (Lemma \ref{algebraisch}), claim ($\ast$) follows 
   with $\Unull := \U$ and $G_0 := \delta_{\nu} f_{\nu}$.

   Let now $\tau_0 \in \Ha$. 
   By assumption we then have 
   $\tau_0 \in \Ha - ({\cal S} \cup \Seins) = \Ha - \Seins$, and  
   Lemma \ref{salpha} implies that  
   $\langle \tau_0, \alpha^{\vee} \rangle = t_0 \in \mathbb Z$ for exactly
   one $\alpha \in R^+$ and  
   $\mu_0 =\nu - \langle \nu +\rho +\tau_0,\alpha^{\vee}\rangle\alpha$ is
   a weight of $E$ with $\mu_0 > \nu$.

   Again we construct a set $\Uhmu$ with a particular choice of 
      $\{h_{\mu}\}$. Namely, choose for all weights
      $\mu \in P(E)$ with $\mu > \nu$ and $\mu \neq \mu_0$ elements
      $h_{\mu} \in \f{h}$, such that  
      $H_{\mu}(w \cdot (\tau_0 +\nu),\tau_0) = 
      \langle w \cdot (\tau_0 +\nu) -\tau_0 -\mu,h_{\mu}\rangle\neq 0$ for all
      $w \in \cal W$. 
   For $\mu_0$ however, choose $h_{\mu_0}$ such that
    \begin{itemize}
    \item[a)] $\langle \alpha, h_{\mu_0} \rangle \neq 0$ and
    \item[b)] $H_{\mu_0}(w \cdot (\tau_0 +\nu),\tau_0) =
              \langle w \cdot (\tau_0 +\nu) -\tau_0 -\mu_0,h_{\mu_0}\rangle
              \neq 0$ for all $s_{\alpha} \neq w \in \cal W$.
    \end{itemize}
   Lemma \ref{salpha} ensures that such a choice of $\{ h_{\mu}\}$ is always
      possible.
   For $\tau \in \f{h}^*$ we set again 
      $p_{\tau} := \prod_{\mu \in P(E), \mu > \nu} \sym H_{\mu}(-,\tau)$, 
      and 
      ${\xi}^{-1}(p_{\tau}) =:z_{\tau} \in \f{Z}$.
   For $u$ defined by
      $u(\tau) = p_{\tau}(\tau +\nu)$ we then obtain
      $\Uhmu = \{ \tau \in \f{h}^* \mid  u(\tau) \neq 0 \}$.
   According to Lemma \ref{resultat} the map $f_{\nu}$ restricted to
      $\Uhmu$ is algebraic.
   Let's have a closer look at the map $u$ :
   It vanishes along the hyperplane $\ker H_{\alpha,t_{0}}$, because 
      for all $\tau \in \f{h}^*$ we have
      ${H_{\mu_0}\big(s_{\alpha} \cdot(\tau +\nu),\tau\big)}
        = -\langle \alpha,h_{\mu_0} \rangle \cdot 
            \big( \langle \tau,\alpha^{\vee} \rangle -t_{0} \big)
        = a \cdot H_{\alpha,t_{0}}(\tau)$,
      where $a := -\langle \alpha,h_{\mu_0} \rangle \neq 0$ according to
      assumption a).

   All the other hyperplanes along which $u$ vanishes are not equal
      to $\ker H_{\alpha,t_{0}}$ but do at most intersect it. 
   Otherwise (since $\tau_0 \in \ker H_{\alpha,t_{0}}$) we would have
      $H_{\mu}\big( w \cdot(\tau_0 +\nu),\tau_0\big)  = 0$ for a pair
      $(\mu,w) \neq (\mu_0, s_{\alpha})$ which is impossible by our
      choice of the $h_{\mu}$.

   For $\bar u: \, \f{h}^* \to k$ defined by
      $\bar u H_{\alpha,t_{0}} = u$, it follows by what we just said that
      $\bar u(\tau_0) \neq 0$.
   With $\Unull := \{ \tau \in \f{h}^* \; |\; \bar u(\tau) \neq 0\}$ 
      we then get
 \begin{itemize}
    \item $\tau_0 \in \Unull$
    \item $1/ {\bar u}$ is algebraic on $\Unull$
    \item $\Uhmu = \Unull - \ker H_{\alpha, t_0}$, thus
          $\Unull = \Uhmu \cup \ker H_{\alpha, t_0}$
    \item $(\U \cap \Unull) \subset \Uhmu$, since 
          $\U \cap \ker H_{\alpha, t_0} = \emptyset$.
 \end{itemize}

   If we define 
      $\bar \delta =\bar \delta_{\nu}:\, \f{h}^* \to k$ by
      $\bar \delta H_{\alpha,t_{0}} = \delta_{\nu}$ and set
             \[
             \begin{array}{lccc}
              G_0 :\, & \Unull & \longrightarrow & V\\
                      & \tau    & \mapsto & \bar \delta(\tau)  
                                  \big(\bar u(\tau)\big)^{-1} 
                                  z_{\tau} (e_{\nu} \otimes v_{\tau}) 
             \end{array}
             \]            
      we know that $G_0$ is algebraic on $\Unull$. 
   In particular, $G_0$ is algebraic on an open neighborhood of $\tau_0$. 
   This is now the local algebraic extension of $\delta_{\nu} f_{\nu}$
      around $\tau_0$ which we were looking for, because according to 
      Lemma \ref{resultat} we have for all $\tau \in \Uhmu$:  
   \begin{eqnarray*}
    \delta_{\nu}(\tau) f_{\nu}(\tau) 
    &=& \delta_{\nu}(\tau)\left(u(\tau)\right)^{-1} 
                               z_{\tau} (e_{\nu} \otimes v_{\tau})\\
    &=& \bar \delta_{\nu}(\tau) H_{\alpha,t_0}(\tau)
          \left( H_{\alpha,t_0}(\tau)\right)^{-1}
          \left(\bar u(\tau) \right)^{-1}
         z_{\tau} (e_{\nu} \otimes v_{\tau})\\
    &=& G_0(\tau).
   \end{eqnarray*}
    In particular, this equation holds for all
       $\tau \in (\U \cap \Unull) \subset \Uhmu$, i.e. ($\ast$).\\ 


  It remains to show that the set of zeros of $G$ has codimension $\ge 2$.

  In order to do this, let $\Szwei$ be the union of all intersections of 
     hyperplanes in $\Ha$ with any other integral hyperplanes, that is
     \[\Szwei := \Seins \cup 
              \{ \tau \in \Ha \mid \langle \tau, \alpha^{\vee} \rangle =m \in
                                \mathbb Z \mbox{ for a pair } (\alpha, m) 
                                \notin \N \} .\]

  We claim
     \begin{itemize}
     \item[($\ast \ast$)] ${\cal K} := \{\tau\in\f{h}^*\;|\; G(\tau)=0 \}
                        \subset (\cal S \cup \Szwei )$
     \end{itemize}
     and conclude then that $\codim \cal K \ge 2$. 
  Indeed, by definition the set $\cal S \cup \Szwei$ consists of a 
     countable family of intersections of integral hyperplanes. 
  On the other hand, $G$ is algebraic on $\f{h}^*$ and therefore
     its set of zeros must be Zariski closed.
  Since a countable family of intersections of integral hyperplanes is
     Zariski closed only if it is a {\em finite} family of such intersections,
     the set $\cal K$ must have codimension $\ge 2$. 
 
  Let now $\tau_0 \in \f{h}^* - (\cal S \cup \Szwei )$. 
  We have to show $G(\tau_0) \neq 0$.

  If $\tau_0 \notin \Ha$, we have $\delta_{\nu}(\tau_0) \neq 0$ and
     $\tau_0 \in \U = \f{h}^*-(\Ha \cup {\cal S})$.
  Now $G$ equals $\delta_{\nu} f_{\nu}$ on $\U$ and
     hence $G(\tau_0) = \delta_{\nu}(\tau_0) f_{\nu}(\tau_0) \neq 0$.

  Let now $\tau_0 \in \Ha$.
  Again by the first part of the proof we know that there is an open 
     neighborhood ${\cal U}_0$ of $\tau_0$, such that  
  $G(\tau) = \bar \delta_{\nu} (\tau) \big(\bar u(\tau)\big)^{-1} 
                                z_{\tau} (e_{\nu} \otimes v_{\tau})$
  for all $\tau \in \Unull$.
  Since $\bar \delta_{\nu} (\tau_0) \neq 0$ it thus suffices to show 
  \begin{itemize}
     \item[($\ast\ast$)']
        $z_{\tau_0} (e_{\nu} \otimes v_{\tau_0}) \neq 0$.
  \end{itemize}

  In order to do this, we will use Lemma \ref{PO} to
     \ref{Ann} of the following section.
  Set $e:= e_{\nu}$ and denote the projection onto the central character 
     $\chi(\tau_0 +\nu)$ by $\pr = \pr_{\chi(\tau_0 +\nu)}:\, V 
     \twoheadrightarrow V$. 

  Since $\tau_0 \in \Ha - \Szwei$, there exists exactly one $\alpha \in R^+$
     such that $\langle \tau_0, \alpha^{\vee} \rangle \in \mathbb Z$. 
  Therefore, 
    $R_{\tau_0} 
    :=\{ \beta \in R \mid \langle \tau_0,\beta^{\vee}\rangle \in \mathbb Z\}
     = \{\alpha, -\alpha \}$ 
     and hence the Weyl group ${\cal W}_{\tau_0}$ of this root system 
     $R_{\tau_0}$ equals 
     ${\cal W}_{\tau_0}
     = \langle s_{\beta} \mid \beta \in R_{\tau_0} \rangle
     = \langle s_{\alpha} \rangle$.
  For $\tau_0$ and 
     $\mu_0 = \nu - \langle \tau_0 +\nu +\rho, \alpha^{\vee} \rangle \alpha$
     we may then apply Lemma \ref{KetteR} and obtain a short exact sequence
     \[{\bigoplus}_{\dim E_{\mu_0}} \verma{s_{\alpha} \cdot (\tau_0 +\nu)}
     \hookrightarrow \pr \big(E \otimes \verma{\tau_0}\big) 
     \twoheadrightarrow L(\tau_0 +\nu)\] 
     such that $\pr(e \otimes v_{\tau_0})/\pr \big
                                            (E \otimes \verma{\tau_0}\big)$
     generates the simple module $L(\tau_0 +\nu)$.
  The module $\pr\big(E \otimes \verma{\tau_0}\big)$ is projective
     (Lemma \ref{verma-einfach}) and this forces the above short exact 
     sequence to be a nontrivial extension because otherwise, being a direct 
     summand, $L(\tau_0 +\nu)$ would be projective too. 
  This contradicts Lemma \ref{verma-einfach}.
 
  Therefore (see for example \cite[Ch.3, Lemma 4.1]{hs}) we have
     at least one direct summand $\verma{s_{\alpha} \cdot (\tau_0 +\nu)}$ of 
     $R_2:= 
       {\bigoplus}_{\dim E_{\mu_0}} \verma{s_{\alpha} \cdot (\tau_0 +\nu)}$,
     such that the projection 
     $\p: R_2 \twoheadrightarrow \verma{s_{\alpha} \cdot (\tau_0 +\nu)}$
     induces a nontrivial extension
     \[ \verma{s_{\alpha} \cdot (\tau_0 +\nu)} \hookrightarrow
          \PO \twoheadrightarrow L(\tau_0 +\nu)\]
     together with a homomorphism
     $\tilde {\p} : R_1:=
                   \pr\big(E\otimes\verma{\tau_0}\big) \to \PO$,
     such that the following diagram commutes
     (the rows are exact, $\PO$ is a push-out of  
     $\p: R_2 \twoheadrightarrow \verma{s_{\alpha} \cdot (\tau_0 +\nu)}$ 
     and $R_2 \hookrightarrow R_1$, see for example the dual 
     statement to \cite[Ch.3, Lemma 1.3]{hs}):
          \[
          \begin{CD}
          R_2        @>>>     R_1         @>>>   R_1/R_2 \cong L(\tau_0 +\nu)\\
          @VV{\p}V            @VV{\tilde {\p}}V    @VV{\id}V\\
          \verma{s_{\alpha} \cdot (\tau_0 +\nu)} @>>> \PO @>>> 
                \PO/\verma{s_{\alpha} \cdot (\tau_0 +\nu)} \cong L(\tau_0 +\nu)
          \end{CD}
          \]
  Since the bottom row is a nontrivial extension and since 
     ${\cal W}_{\tau_0 +\nu} = {\cal W}_{\tau_0} = \langle s_{\alpha} \rangle$,
     the push-out $\PO$ is isomorphic to $P(\tau_0 +\nu)$, the projective 
     cover of $L(\tau_0 +\nu)$ in $\cal O$ (Lemma \ref{PO}).
  In the right loop of the above diagram lies the following 
     element diagram
      \[\xymatrix{
         {\pr(e \otimes v_{\tau_0})} \ar@{|->}[d] \ar@{|->}[r]
        &{\pr(e \otimes v_{\tau_0})/R_2} \ar@{|->}[d]\\
         {\tilde {\p} \big(\pr(e \otimes v_{\tau_0})\big)}\ar@{|->}[r]
        &{\tilde {\p} \big(\pr(e \otimes v_{\tau_0})\big)/
                   \verma{s_{\alpha} \cdot (\tau_0 +\nu)}}
        }\]
   Since $\pr(e \otimes v_{\tau_0})/R_2$ is a generator of $L(\tau_0 +\nu)$
      we conclude that
      $\tilde {\p} \big(\pr(e \otimes v_{\tau_0})\big)/
                                    \verma{s_{\alpha} \cdot (\tau_0 +\nu)} 
      \in \PO/\verma{s_{\alpha} \cdot (\tau_0 +\nu)}$
      is also a generator of $L(\tau_0 +\nu)$ and
      thus its preimage $\tilde {\p} \big(\pr(e \otimes v_{\tau_0})\big)$ 
      generates the indecomposable module $\PO \cong P(\tau_0 +\nu)$.
   By Lemma \ref{Ann} the element $z_{\tau_0}$ is not contained in the
      annihilator of $P(\tau_0 +\nu)$, this means in particular  
      $0 \neq z_{\tau_0} \tilde {\p} \big(\pr(e \otimes v_{\tau_0})\big)
         = \tilde {\p} \big(z_{\tau_0} \pr(e \otimes v_{\tau_0})\big)$
      and we obtain $z_{\tau_0} \pr(e \otimes v_{\tau_0}) \neq 0$. 
   Therefore, claim ($\ast\ast$)' is proven and the proof of 
      Theorem \ref{theoeins} is finished provided we know that Lemma \ref{PO},
      \ref{verma-einfach},  \ref{KetteR} and \ref{Ann} hold.
\end{proof}

\subsection{Proof of Lemma \ref{PO}, \ref{verma-einfach}, \ref{KetteR} 
             and \ref{Ann}}\label{lemmata}

In the following let $e_{\nu}$ denote again the fixed extremal weight
  vector of the finite dimensional irreducible $\f{g}$-module $E = E(\nu)$.
Set $v_{\tau} \in \verma{\tau}$ the canonical generator of the Verma module
  and let $V := E \otimes \f{U}(\f{n}^-) \cong 
  E \otimes \verma{\tau}$ for all $\tau \in \f{h}^*$.
The category $\cal O$ is the category of all finitely generated 
  $\f{g}$-modules, which are $\f{b}$-finite and semisimple over $\f{h}$ 
  \cite{bgg2}.
For $\lambda \in \f{h}^*$ set 
   ${\cal W}_{\lambda} := \{w \in {\cal W} \mid \lambda - w \lambda \in P\}$
   and $R_{\lambda} 
   :=\{ \beta \in R \mid \langle \lambda,\beta^{\vee}\rangle \in \mathbb Z\}$.
The group ${\cal W}_{\lambda}$ is then the Weyl group to the root system
   $R_{\lambda}$ \cite[1.3]{ja1}.
Denote by $P(\lambda)$ the projective cover in $\cal O$ of the simple 
   module $L(\lambda)$.

\begin{samepage}
\begin{lemma}\label{PO}
      Let $\lambda \in \f{h}^*$ with 
      ${\cal W}_{\lambda} = \langle s_{\alpha} \rangle$ and
      $\verma{\lambda}=L(\lambda)$ simple.
      If the short exact sequence 
      $\verma{s_{\alpha} \cdot \lambda} \hookrightarrow N 
         \twoheadrightarrow L(\lambda)$
      is a nontrivial extension, 
      then $N$ is isomorphic to $P(\lambda)$. 
\end{lemma} 
\end{samepage}
%
\begin{proof}
      For $\lambda$ with 
      ${\cal W}_{\lambda} = \langle s_{\alpha} \rangle$ there exists in 
      $\cal O$ (up to isomorphism) a unique indecomposable projective
      module $P(\lambda)$, such that the short exact sequence
      \[\verma{s_{\alpha} \cdot \lambda} \hookrightarrow P(\lambda) 
      \twoheadrightarrow L(\lambda)\]
      is a nontrivial extension \cite[Ch.4]{bgg2}. 
      We thus have a nontrivial element of 
      $\Ext^1\big(L(\lambda),\verma{s_{\alpha} \cdot \lambda}\big)$.
      The assertion follows immediately if we knew that 
      $\dim \Ext^1\big(L(\lambda),\verma{s_{\alpha}\cdot \lambda}\big) \le 1$.
      Indeed, construct for 
      $ \verma{s_{\alpha} \cdot \lambda} \hookrightarrow P(\lambda) 
      \twoheadrightarrow L(\lambda)$
      the long exact homology sequence \cite[Ch.3, Thm.5.3]{hs}
      (set $M=\verma{s_{\alpha} \cdot \lambda}, P = P(\lambda), 
      L = L(\lambda), \Ext = \Ext_{\cal O}$ and let
      $\omega$ be the connecting homomorphism):
      \[ \cdots \to \Hom(M,M) \stackrel{\omega}{\to}
         \Ext^1(L,M) \to \Ext^1(P,M) \to \Ext^1(M,M) \to \cdots \]
      and note that $\Ext^1(P,M)=0$ since $P = P(\lambda)$ is projective
      \cite[Chap.3, Prop.2.6.]{hs}. 
      By exactness we obtain $\omega$ surjective.
      Since now for any Verma module M $\dim \Hom(M,M) = 1$, 
      we conclude that $\dim \Ext^1(L,M) \le 1$. 
\end{proof} 

\begin{samepage}
\begin{lemma} \label{verma-einfach}
         Let $\tau \in \f{h}^*$ with
         ${\cal W}_{\tau} = \langle s_{\alpha} \rangle$.
\begin{itemize}
 \item[(a)] If $\langle \tau +\rho, \alpha^{\vee} \rangle \ge 0$, then
            $\pr_{\chi(\tau +\nu)}\big(E \otimes \verma{\tau}\big)$\;
            is projective, this means
            $\pr_{\chi(\tau +\nu)}\big(E \otimes \verma{\tau}\big)$ 
            is a projective object of the category $\cal O$.
 \item[(b)] If $\langle \tau +\nu +\rho, \alpha^{\vee} \rangle \le 0$,
            then $\verma{\tau +\nu} = L(\tau +\nu)$ is simple and not
            projective.
\end{itemize}
\end{lemma}
\end{samepage}
\begin{proof}
         a) A Verma module $\verma{\lambda}$ is projective if
            $\langle \lambda +\rho, \beta^{\vee} \rangle \ge 0$ for all
            $\beta \in R_{\lambda} \cap R^+$ \cite[4.8]{ja2}. 
            Since ${\cal W}_{\tau} = \langle s_{\alpha} \rangle$,
            we have $R_{\tau} \cap R^+ = \{ \alpha \}$ and hence 
            \verma{\tau} is projective.
            Then also $E \otimes \verma{\tau}$ is projective, 
            because for $E$ with $\dim E < \infty$ the functor
          \begin{center}
             $ 
             \begin{array}{lccc}
              F_E :\, & \cal O & \longrightarrow & \cal O\\
                      & M      & \mapsto         & E \otimes M
             \end{array}
             $
         \end{center}

            maps projective objects of $\cal O$ to projective objects of
            $\cal O$ \cite{bg}.
            Its direct summand  
            $\pr_{\chi(\tau +\nu)}\big(E \otimes \verma{\tau}\big) \subset
            E \otimes \verma{\tau}$ is then projective too.

         b) A Verma module $\verma{\lambda}$ is simple if and only if
            $\langle \lambda + \rho, \beta^{\vee} \rangle \le 0$ for all
            $\beta \in R_{\lambda} \cap R^+$ \cite[7.6.24]{di} or
            \cite[1.8, 1.9]{ja1}.
            Since ${\cal W}_{\tau +\nu} = {\cal W}_{\tau} = 
            \langle s_{\alpha} \rangle$ we obtain
            $R_{\tau+\nu} \cap R^+ = \{\alpha\}$ and the simplicity of  
            $\verma{\tau +\nu} = L(\tau +\nu)$ follows.
            Since the short exact sequence 
            $\verma{s_{\alpha} \cdot (\tau +\nu)} \hookrightarrow 
            P(\tau +\nu) \twoheadrightarrow L(\tau +\nu)$ 
            is a nontrivial extension \cite[Ch.4]{bgg2}, 
            the module $L(\tau +\nu)$ cannot be projective.
\end{proof}

\begin{samepage}
\begin{lemma}\label{KetteR}
      Let $\tau \in \f{h}^*$ with ${\cal W}_{\tau} =\langle s_{\alpha} \rangle$
      and
      $\mu_0 := 
               \nu - \langle \tau +\nu +\rho, \alpha^{\vee} \rangle \alpha
               \in P(E)$ such that $\mu_0 > \nu$.
      Then the module 
      $\pr_{\chi(\tau +\nu)}\big(E \otimes \verma{\tau}\big)$ 
      admits a chain of submodules 
      \[\pr_{\chi(\tau +\nu)}\big(E \otimes \verma{\tau}\big)=R_1 
                          \supset R_2 \supset 0,\]
      such that $R_2 \cong {\bigoplus}_{\dim E_{\mu_0}} 
                \verma{s_{\alpha} \cdot (\tau +\nu)}$ and
      $R_1/R_2 \cong \verma{\tau +\nu} = L(\tau +\nu)$ and
      such that $L(\tau +\nu)$ is generated by 
      $\big(\pr_{\chi(\tau +\nu)}(e_{\nu} \otimes v_{\tau})\big)/R_2$.
\end{lemma} 
\end{samepage}

\begin{proof}
  Let $\mu_1,\mu_2,\ldots,\mu_m$ be the weights of $E$ without multiplicities,
  numbered such that $\mu_i < \mu_j \Rightarrow i < j$. 
  Set $\mu_{j_0} := \nu, \mu_{j_1} := \mu_0$, $d_j:=\dim E_{\mu_j}$ and 
  $\pr := \pr_{\chi(\tau +\nu)}$.
  From chapter 2 we already know that $E \otimes \verma{\tau}$ has a chain of
  submodules
  $E \otimes \verma{\tau} =M_1 \supset \cdots \supset M_m \supset 0$
  such that 
  $M_j /M_{j+1}\cong {\bigoplus}_{d_j} \verma{\tau +\mu_j}$.
  It is then clear that 
  $\pr(E \otimes \verma{\tau}) =\pr M_1 \supset \cdots \supset 
                                                          \pr M_m \supset 0$
  is a chain of submodules such that the subquotients 
  $\pr M_j /\pr M_{j+1}$ are isomorphic to
  $\bigoplus_{d_j} \verma{\tau +\mu_j}$,
  if $\tau +\mu_j \in {\cal W} \cdot (\tau +\nu)$, otherwise they are $0$. 
  Since $\mu_j$ and $\nu$ are integral weights,
  we get $\mu_j = w \cdot (\tau + \nu) -\tau$ only for $w \in \cal W_{\tau}$.
  For $w = e$ this yields $\mu_{j_0} = \nu$, for $w = s_{\alpha}$ we obtain
  $\mu_{j_1} = \mu_0$ since $\tau +\mu_0 = s_{\alpha} \cdot (\tau +\nu)$. 
  By assumption we have $\mu_{j_0} = \nu < \mu_0 = \mu_{j_1}$, 
  and thus $j_0 < j_1$.
  Therefore, if we omit in this chain trivial submodules we obtain a chain 
  $\pr\big(E \otimes \verma{\tau} \big) = \pr M_{j_0} =: R_1 
                                     \supset \pr M_{j_1} =: R_2 \supset 0$
  such that
  $R_2 \cong {\bigoplus}_{d_{j_1}} \verma{\tau +\mu_0} 
               = {\bigoplus}_{d_{j_1}} \verma{s_{\alpha} \cdot (\tau +\nu)}$
  and $R_1/R_2 \cong {\bigoplus}_{d_{j_0}}
                \verma{\tau +\nu} \cong \verma{\tau +\nu}$.
  By construction of the $M_i$ this module is generated by 
  $\pr_{\chi(\tau +\nu)}(e_{\nu} \otimes v_{\tau})/R_2$.
  Now $\mu_0 > \nu$ forces 
  $\langle \tau +\nu +\rho, \alpha^{\vee} \rangle < 0$ and together with 
  ${\cal W}_{\tau +\nu} = {\cal W}_{\tau} =\langle s_{\alpha} \rangle$ 
  this implies $\verma{\tau +\nu} = L(\tau +\nu)$ (Lemma \ref{verma-einfach}).
\end{proof}

For $\mu \in \f{h}^*$ define maps
   $H_{\mu}: \, \f{h}^* \times \f{h}^* \to k$ by
   $H_{\mu}(\lambda,\tau) := \langle \lambda - \tau -\mu,h_{\mu} \rangle
                   \; \forall \lambda, \tau \in \f{h}^*, h_{\mu} \in \f{h}$.
For $\tau \in \f{h}^*$ set 
   $p_{\tau} := \prod_{\mu > \nu, \mu \in P(E)} \sym H_{\mu}(-,\tau)$
   and $z_{\tau} :={\xi}^{-1}(p_{\tau}) \in \f{Z}$ the preimage of $p_{\tau}$
   under the Harish-Chandra isomorphism  
   $\xi: \,\f{Z} \stackrel{\sim}{\longrightarrow}
                        {\cal P}(\f{h}^*)^{\cal W \cdot}$.
\begin{samepage}
\begin{lemma} \label{Ann}
     Let $\tau_0 \in \f{h}^*$ with
         ${\cal W}_{\tau_0} = \langle s_{\alpha} \rangle$ 
     and
         $-\langle \rho,\alpha^{\vee} \rangle \le 
                            \langle \tau_0,\alpha^{\vee} \rangle \le 
                           -\langle \nu+\rho, \alpha^{\vee} \rangle - 1$.
     Assume the $h_{\mu} \in \f{h}$ to be chosen such that
     \begin{itemize}
           \item[(a)] $\langle \alpha, h_{\mu_0} \rangle \neq 0$ for 
                   $\mu_0 := 
                   \nu - \langle \tau_0+\nu+\rho,\alpha^{\vee} \rangle \alpha$
            and
            \item[(b)] $H_{\mu}\big(w \cdot(\tau_0+\nu),\tau_0\big) \neq 0$ 
                    for $(w,\mu) \neq (s_{\alpha},\mu_0), 
                    \forall w \in {\cal W}, \mu \in P(E)$ with $\mu > \nu$.
     \end{itemize} 
     Then 
     $z_{\tau_0} \notin \Ann_{\f{Z}}P(\tau_0 +\nu)$.
\end{lemma} 
\end{samepage}
%
\begin{proof}
   We will first give an explicit description of the annihilator 
      $\Ann_{\f{Z}}P(\tau_0 +\nu)$ by a theorem of Soergel \cite{so}.
   For this let $\lambda \in \f{h}^*$, such that for all
      $\beta \in R^+\cap R_{\lambda}$ we have  
      $\langle \lambda +\rho,\beta^{\vee} \rangle \ge 0$.
   Denote by 
      $w_{\lambda} \in \cal W_{\lambda}$ the longest element with respect to
      the Bruhat ordering.
   Then define for $\mu \in \f{h}^*$ the map 
             \[
             \begin{array}{lccc}
              \mu^+:\, &\cal P(\f{h}^*) & \longrightarrow & \cal P(\f{h}^*)\\
                       & p              & \mapsto         & \mu^+(p)
             \end{array}
             \]
       by $\big(\mu^+(p)\big)(\tau) := p(\tau +\mu)$\; for all 
       $\tau \in \f{h}^*$.
   Then \cite[2.2]{so}:
       \[
          {\xi}^{-1}(p) \in \Ann_{\f{Z}}P(w_{\lambda} \cdot \lambda) 
          \Longleftrightarrow
          \lambda^+(p) \in 
          \big(\cal P^+(\f{h}^*)\big)^{\cal W_{\lambda}} \cal P(\f{h}^*).
       \]
   Here $\cal P^+(\f{h}^*)$ denotes the set of polynomial functions on
       $\f{h}^*$ without constant term.
   We want to describe the annihilator of $P(\tau_0 +\nu)$ and choose 
      for this $\lambda_0 := s_{\alpha} \cdot(\tau_0 +\nu)
      = \tau_0 +\mu_0
      = \tau_0 +\nu - \langle \tau_0 +\nu +\rho, \alpha^{\vee}\rangle \alpha$.
   We then have
      ${\cal W}_{\lambda_0}={\cal W}_{s_{\alpha} \cdot(\tau_0 +\nu)} = 
      {\cal W}_{\tau_0 +\nu} = {\cal W}_{\tau_0} = \langle s_{\alpha} \rangle$
      and hence $R_{\lambda_0} = \{ \alpha, -\alpha \}$.
   By assumption we know that
      $\langle \lambda_0 +\rho, \alpha^{\vee} \rangle \ge 0$ and obtain
      for all $\beta \in R^+ \cap R_{\lambda_0}$ that 
      $\langle \lambda_0 +\rho,\beta^{\vee} \rangle \ge 0$. 
   Now $w_{\lambda_0} = s_{\alpha}$ is the longest element in 
      ${\cal W}_{\lambda_0}$ and we may apply Soergel's theorem with 
      $\lambda:= \lambda_0= s_{\alpha} \cdot (\tau_0 +\nu)$.
   Assume we had 
      ${\xi}^{-1}(p)\in \Ann_{\f{Z}}P(w_{\lambda_0} \cdot\lambda_0)
                                            = \Ann_{\f{Z}}P(\tau_0 +\nu)$ 
      for $p = \sum_{i=1}^{n} p_i q_i$ with 
      $p_i \in \cal P^+(\f{h}^*)\big)^{\cal W_{\lambda_0}}$ and
      $q_i \in \cal P(\f{h}^*)$.
   Then it follows for all $p_i$ that 
      \[p_i(\lambda_0 + \mu) 
              = \left(\lambda_0^+(p_i)\right)(\mu) 
              =\left(\lambda_0^+(p_i)\right)(s_{\alpha} \mu) 
              = p_i(\lambda_0 + s_{\alpha}\mu) \; \forall \mu \in \f{h}^*.\]
   For $\mu = \alpha$ we obtain 
      $p_i(\lambda_0 + \alpha) = p_i(\lambda_0 - \alpha)$ and
      this forces the derivative of $p$ in direction of $\alpha$ to vanish
      at the point $\lambda_0$. 
   Let's check this condition for $p_{\tau}$.     
   By definition we have for all $\lambda \in \f{h}^*$:
      $p_{\tau}(\lambda)
       =\prod_{\mu > \nu , w \in \cal W} H_{\mu}(w \cdot \lambda,\tau)$.
   If we define $\bar p_{\tau}$ by
      $\bar p_{\tau} H_{\mu_0}(-,\tau) = p_{\tau}$, we obtain
      \[ \bar p_{\tau}(\lambda) := 
          \prod_{\substack{(w,\mu) \neq (e, \mu_0)\\
                           w \in {\cal W}, \mu > \nu}}
                                   H_{\mu}(w \cdot \lambda,\tau).\]
   Let $p_{\tau}^{\prime}$ denote the derivative of $p_{\tau}$ in
      direction $\alpha$. By the product rule we have
      $p_{\tau}^{\prime} =
      \bar p_{\tau}^{\prime} H_{\mu_0}(-,\tau) + 
                          \bar p_{\tau} H_{\mu_0}^{\prime}(-,\tau)$.
     
   Since $s_{\alpha} \cdot (\tau_0 +\nu) = \tau_0 +\mu_0$ it follows that
      $H_{\mu_0}\big(s_{\alpha} \cdot (\tau_0 +\nu),\tau_0\big) =0$
      and we get
      $p_{\tau_0}^{\prime}\big(s_{\alpha} \cdot (\tau_0 +\nu)\big) =
        \bar p_{\tau_0}\big(s_{\alpha} \cdot (\tau_0 +\nu)\big) 
        H_{\mu_0}^{\prime}\big(s_{\alpha} \cdot (\tau_0 +\nu),\tau_0\big)$.
   But now neither of these two factors is zero because
      \begin{eqnarray*}
           \bar p_{\tau_0}\big(s_{\alpha} \cdot (\tau_0 +\nu)\big)
           &=&
           \prod_{\substack{(w,\mu) \neq (e, \mu_0)\\
                                         w \in {\cal W}, \mu > \nu}}
             H_{\mu}\big(w \cdot s_{\alpha} \cdot (\tau_0 +\nu),\tau_0\big)\\
           &=&
           \prod_{\substack{(w,\mu) \neq (s_{\alpha}, \mu_0)\\
                                         w \in {\cal W}, \mu > \nu}}
                 H_{\mu}\big(w \cdot (\tau_0 +\nu),\tau_0\big)
      \end{eqnarray*}
      and by assumption (b) none of the factors 
      $H_{\mu}\big(w \cdot (\tau_0 +\nu),\tau_0\big)$ vanishes, hence 
      $\bar p_{\tau_0}\big(s_{\alpha} \cdot (\tau_0 +\nu)\big) \neq 0$.
   On the other hand, also 
      $H_{\mu_0}^{\prime}\big(s_{\alpha}\cdot (\tau_0 +\nu),\tau_0\big)\neq 0$
      since by assumption (a) we have:
      $H_{\mu_0}^{\prime}(\lambda,\tau_0) 
           = \ein{\frac{\dif}{\dif t}}{t=0} 
               \langle \lambda +t \alpha - \tau_0 - \mu_0,h_{\mu_0} \rangle
           = \langle \alpha,h_{\mu_0} \rangle
           \neq 0$.

  Together this implies
     $p_{\tau_0}^{\prime}\big(s_{\alpha} \cdot (\tau_0 +\nu)\big) \neq 0$
     and therefore $\xi(p_{\tau_0}) = z_{\tau_0}$ cannot be contained 
     in the annihilator of $P(\tau_0 +\nu)$.
\end{proof}


     \section{\bf The triangle function $\boldsymbol \Delta$} \label{kapvier}

\subsection{Preliminaries} 

  Let $M$ be a representation of $\f{g}$ and $E$ a vector space. 
  Then $E \otimes M$ is a representation of $\f{g}$ via
    $X(e \otimes m) := e \otimes Xm$ for all $X \in \f{g}, e \in E$ and 
    $m \in M$.
  If in addition $E$ is also a representation of $\f{g}$, then we obtain
    a second $\f{g}$-operation on $E \otimes M$ via 
    $X(e \otimes m) := Xe \otimes m + e \otimes Xm$. 
  To distinguish these two representations we denote the first one by
    $E \hotimes M$.
  Let now $E$ be a finite dimensional representation of $\f{g}$ and 
    $\nu \in P$ a weight of $E$. 

\begin{lemma} \label{can}
  Let $\tau \in \f{h}^*$ such that $\chi(\tau +\nu) \neq 
   \chi(\tau +\mu)$ for all $\mu \in P(E)$ with $\mu \neq \nu$.
  Then there exists a unique natural isomorphism
  \[\can: \; E_{\nu} \hotimes \verma{\tau +\nu} 
              \stackrel{\sim}{\longrightarrow}
          \pr_{\chi(\tau +\nu)}(E \otimes \verma{\tau}) ,\]
    such that $\can (e \hotimes v_{\tau +\nu}) =
         \pr_{\chi(\tau +\nu)}(e \otimes v_{\tau})$ 
         for all $e \in E_{\nu}$.
\end{lemma}

\begin{bemerkung}
  For a {\em generic} weight, i.e. a weight $\tau \in \f{h}^*$ such that
    $\langle \tau, \alpha^{\vee} \rangle \notin \mathbb Z$ for all
    $\alpha \in R$, the central characters $\chi(\tau +\mu)$ for 
    $\mu \in P(E)$ are pairwise distinct.
  In particular, in this case the condition of the lemma is always satisfied.
\end{bemerkung}

\begin{proof}
  By the so-called tensor identity we have a canonical isomorphism
   \[ \f{U}\otimes_{{\f{U}}(\f{b})}(E \otimes k_{\tau}) 
      \stackrel{\sim}{\longrightarrow}
       E \otimes (\f{U} \otimes_{{\f{U}}(\f{b})} k_{\tau}) \]
   such that 
   $u \otimes (e \otimes a) \mapsto u(e \otimes(1 \otimes a))$.
  Call the left hand side $F$, the right hand side is $E \otimes \verma{\tau}$.
  Denote by $\mu_1, \ldots, \mu_m$ the weights of $E$. 
  The filtration
    $E \otimes \verma{\tau} =M_1 \supset \cdots \supset M_m \supset 0$, 
    where the subquotients are isomorphic to direct sums of Verma modules
    (see page \pageref{filtration}), induces a filtration of $F$ : 
    $\f{U}\otimes_{{\f{U}}(\f{b})}(E \otimes k_{\tau}) = 
    F=F_1 \supset F_2 \supset \cdots \supset F_m\supset 0$ such that
    \[F_j/F_{j+1} \cong 
        \f{U}\otimes_{{\f{U}}(\f{b})}(E_{\mu_j} \otimes k_{\tau}).\]
   Now here the right hand side is canonically isomorphic to  
    $E_{\mu_j} \hotimes \verma{\tau+\mu_j}$ 
    (by mapping $e \otimes u v_{\tau +\mu_j}\mapsto u \otimes (e \otimes 1)$) 
    and in particular, we have that 
    $\chi \left(\tau + \mu_j\right) \left(F_j / F_{j+1} \right) = 0$.

   Let now $\nu=\mu_i$ for a fixed $i$. 
   By the condition on $\tau$ we know that
     $\chi(\tau +\nu) \neq \chi(\tau + \mu_j)$ for all $j \neq i,\: 
     j \in \{1, \ldots ,m\}$ and hence 
     $\pr_{\chi(\tau +\nu)}\left(F_j / F_{j+1} \right) = 0$ for all $j\neq i$.
   Thus we get 
     $\pr_{\chi(\tau +\nu)}F =\pr_{\chi(\tau +\nu)}F_1 = \ldots
                             = \pr_{\chi(\tau +\nu)}F_i$ 
     and also $\pr_{\chi(\tau +\nu)}F_{i+1} = \ldots 
                             = \pr_{\chi(\tau +\nu)}F_m = 0$.
   We conclude that $\pr_{\chi(\tau +\nu)} F \subset F_i$
     and that $F_{i+1} \subset \ker( \pr_{\chi(\tau +\nu)}:   
                             F_i\twoheadrightarrow \pr_{\chi(\tau +\nu)}F)$. 
   Since $\pr_{\chi(\tau +\nu)}\left(F_i / F_{i+1} \right) = F_i / F_{i+1}$, 
     we know even that $F_{i+1}$ is equal to this kernel.
   This now induces a natural isomorphism
     \[ F_i /F_{i+1}\stackrel{\sim}{\longrightarrow} \pr_{\chi(\tau +\nu)}F \]
     such that
     \[
     \begin{array}{ccc}
       E_{\nu} \hotimes \verma{\tau + \nu} & \stackrel{\sim}{\longrightarrow}&
                       \pr_{\chi(\tau + \nu)} F\\
       e \hotimes (u v_{\tau+\nu}) & \longmapsto & 
                       \pr_{\chi(\tau + \nu)} (u \otimes (e \otimes 1))
        \end{array}
     \]

   We apply the tensor identity $F \cong E \otimes \verma{\tau}$ and the
     lemma follows.
\end{proof}

  Let us recall the theorem of Bernstein and Gelfand \cite{bg} for 
    projective functors.
  For this denote by $\cal M$ the category of all $\f{Z}$-finite
    $\f{g}$-modules and for $\chi \in \maxz$ let
    \begin{align*}
    {\cal M}(\chi)&:=
                  \{ M \in {\cal M} \mid \chi M = 0 \} \\ 
    {\cal M}^{\infty} (\chi) &:= \{ M \in {\cal M} \mid \text{ for all }
                    m \in M \text{ exists } n \in \mathbb N 
                    \text{ such that }\chi^n m =0 \}.
    \end{align*} 
   
  A {\em projective $\chi$-functor} is then a functor
    $F:\; {\cal M}(\chi) \to \cal M$, which is isomorphic to a direct summand
    of a functor $E \otimes$ for a finite dimensional representation $E$.
  In particular, the restriction of the translation functor 
     $T_{\tau}^{\tau +\nu}:\,  \cal M^{\infty} \big(\chi(\tau)\big) 
     \to  \cal M^{\infty} \big(\chi(\tau+\nu)\big)$ to the subcategory
     ${\cal M}(\chi(\tau))$ is a projective $\chi(\tau)$-functor.
  For two projective $\chi$-functors 
    $F, \tilde{F}:\; {\cal M}(\chi) \to \cal M$ denote by  
    $\Hom_{{\cal M}(\chi)\to} (F,\tilde{F})$ the space of all natural 
    transformations from $F$ to $\tilde{F}$.

 \begin{theoohne} \cite[3.5]{bg}
     Let $F, \tilde{F} :\; {\cal M}(\chi) \to \cal M$ be projective
       $\chi$-functors and let $\tau \in \f{h}^*$ such that $\chi(\tau)= \chi$ 
       and $\verma{\tau}$ is projective.
     Then the obvious map
       \[\Homtau (F, \tilde{F}) \stackrel{\sim}{\longrightarrow} 
       \Hom_{\f{g}}(F \verma{\tau}, \tilde{F} \verma{\tau})\]
       is an isomorphism.
 \end{theoohne}

 \begin{bemerkung}
   \begin{enumerate}
   \item[(i)]
     The Verma module $\verma{\tau}$ is projective if and only if
      $\langle \tau+ \rho, \alpha^{\vee}\rangle \notin \{-1, -2, \ldots \}$ 
      for all $\alpha \in R^+$. 
     In particular, for a generic weight $\tau$ the Verma module
      $\verma{\tau}$ is always projective.
   \item[(ii)]
     Note that in \cite{bg} the Verma module with highest weight $\tau$ 
      is denoted by $M_{\tau+\rho}$.
     Accordingly, the theorem there is formulated for all $\tau$ with
      $\langle \tau, \alpha^{\vee}\rangle \notin \{-1, -2, \ldots \}$ for 
      all $\alpha \in R^+$.
   \end{enumerate}
 \end{bemerkung}

\subsection{Definition of $\Delta$}

 For $\nu \in P$ let $E(\nu)$ be the finite dimensional irreducible
    $\f{g}$-module with extremal weight $\nu$ and let $x \in \cal W$.
 For a weight $\tau \in \f{h}^*$ with $\chi(\tau +\nu) \neq 
   \chi(\tau +\mu)$ for all $\mu \in P(E)$ with $\mu \neq \nu$, 
   Lemma \ref{can} yields a canonical isomorphism
   \[E(\nu)_{x\nu} \hotimes \verma{x \cdot(\tau+\nu)} 
                   \stackrel{\sim}{\longrightarrow}
   \pr_{\chi(\tau +\nu)}(E(\nu) \otimes \verma{x\cdot\tau}) 
   = T_{\tau}^{\tau +\nu} \verma{x\cdot\tau} .\]
 Let now $\nu, \mu \in P$ be integral weights. 
 We set
   $E^{\prime}:= E(\nu)$, $E^{\prime\prime} := E(\mu)$ and $E := E(\nu +\mu)$. 
 For generic $\tau$ consider the following sequence of 
   isomorphisms: \\

\begin{samepage} 
 \hspace*{0.5em}
  $\Homtau \big(T_{\tau +\nu}^{\tau +\nu +\mu} \circ T_{\tau}^{\tau +\nu},
           T_{\tau}^{\tau +\nu +\mu}\big)$\\[-2.5em]

  \begin{center}
  \begin{align*}
  &\stackrel{\sim}{\longrightarrow} 
     \Hom_{\f{g}}\big(T_{\tau +\nu}^{\tau +\nu +\mu}  T_{\tau}^{\tau +\nu} 
                                                          \verma{x \cdot \tau},
           T_{\tau}^{\tau +\nu +\mu} \verma{x \cdot \tau}\big)\\
  &\stackrel{\sim}{\longrightarrow} 
     \Hom_{\f{g}}\big(E^{\prime\prime}_{x \mu} \hotimes E^{\prime}_{x \nu} 
                                  \hotimes \verma{x \cdot (\tau +\nu +\mu)},
           E_{x(\mu +\nu)} \hotimes \verma{x \cdot (\tau +\nu +\mu)}\big)\\
  &\stackrel{\sim}{\longrightarrow} 
     \Hom_k (E^{\prime\prime}_{x \mu} \otimes E^{\prime}_{x \nu}, E_{x(\mu +\nu)})\\
  &\stackrel{\sim}{\longrightarrow} 
     E^{\prime\prime*}_{x \mu} \otimes E^{\prime*}_{x \nu} \otimes E_{x(\mu +\nu)} 
  \end{align*}
  \end{center}
\end{samepage} 

  Here, we obtain the first isomorphism by the Theorem of Bernstein-Gelfand,
   the second is due to Lemma \ref{can}, the others are obvious.      
  We call this map $\nat(\mu, \nu;x)(\tau)$ \label{nat}
   and define for generic $\tau$ and for the triangle 
   
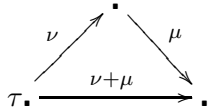
\begin{center}
       {\xymatrix{
        &&&& & {\centerdot} \ar[dr]^{\mu}\\
        &&&&{\tau \centerdot} \ar[ur]^{\nu} \ar[rr]^{\nu+\mu} & 
                                                          & {\centerdot}} }
   \end{center}
   the value of the triangle function $\Delta$ by
   \[\Delta(\mu, \nu;x)(\tau) := \det\left(x^{-1} \circ \nat(\mu, \nu;x)(\tau) 
                               \circ (\nat(\mu, \nu;e)(\tau))^{-1} \right) .\]
 
   We have yet to explain the map
   \[x^{-1}:  E(\mu)^*_{x \mu} \otimes E(\nu)^*_{x \nu} \otimes 
                                              E(\mu+\nu)_{x (\nu +\mu)}
               \to
               E(\mu)^*_{\mu} \otimes E(\nu)^*_{\nu} \otimes 
                                               E(\mu+\nu)_{\nu +\mu} \;.\]

  For this let $G$ be a simply connected algebraic group with Lie algebra 
     $\f{g}$ and $T \subset G$ a maximal torus with Lie algebra $\f{h}$.
  Each finite dimensional representation $E$ of $\f{g}$ is in a natural way
     a representation of $G$. 
  The operation of $N_G (T)$, the normalizer of $T$ in $G$, on $E$ 
     stabilizes $E_0$ and factors over an operation of ${\cal W}= N_G (T)/ T$.
  The map $x^{-1}$ is given by this operation of $\cal W$ on the zero weight
     space $\left(E(\mu)^* \otimes E(\nu)^* \otimes E(\nu +\mu) \right)_0$.
 
  Now take for $\tau$ not only any, but rather {\em the} generic weight:
     For this denote by $S:= S_k(\f{h})$ the symmetric algebra of $\f{h}$
     and let $K:= \Quot(S)$ be its quotient field.
  We then have a $k$-linear map
     $\f{h} \hookrightarrow S \hookrightarrow K$ and obtain thus a 
     $K$-linear map $\tau : K \otimes_k \f{h} \longrightarrow K$
     such that the following diagram commutes 
     \[
     \begin{CD}
        \f{h}               @>>>         S\\
        @VVV                             @VVV\\
        K \otimes_k \f{h}   @>\tau>>     K\\
     \end{CD}
     \]

  With this $\tau$ (= tautologous) we then obtain
          $\Delta (\mu,\nu;x)(\tau) \in K^{\times}$.
  These are the triangle functions.

\subsection{Uncanonical definition of $\Delta$} \label{uncan}

 We defined the triangle functions by a series of canonical isomorphisms.
 For our purposes it is sometimes more convenient to realize the 
   triangle functions in the following -- uncanonical -- way:
 Let $\nu \in P$, $x \in \cal W$ and choose a fixed extremal weight vector 
   $0 \neq e_{\nu} \in E(\nu)_{\nu}$.
 Since extremal weight spaces are one dimensional, this choice is unique 
   up to non-zero scalar.
 By Lemma \ref{can} we obtain for generic $\tau \in \f{h}^*$ 
   a -- no longer canonical -- isomorphism
 \[
  \begin{array}{lccc}
         F_{x\nu}(x\cdot \tau) :& \verma{x \cdot(\tau +\nu)}
                                & \stackrel{\sim}{\longrightarrow} 
                                & T_{\tau}^{\tau +\nu} \verma{x \cdot \tau}\\
                                & v_{x \cdot (\tau +\nu)}
                                & \mapsto
                                & \pr_{\chi(\tau +\nu)}(\dot{x} e_{\nu}
                                             \otimes v_{x \cdot \tau})
  \end{array}
 \]

 Here $\dot{x} \in G$ denotes a pre-image of $x \in {\cal W} \cong N_G T / T$.
 We have $\dot{x} e_{\nu} \in E(\nu)_{x\nu}$.
 Let now $\mu \in P$ be another weight. 
 Choose $\tilde{e}_{\mu} \in E(\mu)_{\mu}$ and 
   $\bar{e}_{\mu +\nu} \in E(\mu +\nu)_{\mu +\nu}$, both non-zero, and 
   consider for generic $\tau$ the following sequence of isomorphisms: \\

 \hspace*{0.5em}
  $\Homtau \big(T_{\tau +\nu}^{\tau +\nu +\mu} \circ T_{\tau}^{\tau +\nu},
                               T_{\tau}^{\tau +\nu +\mu}\big)$ \\[-2.5em]

  \begin{center}
  \begin{align*}
    & \stackrel{\sim}{\longrightarrow} 
      \Hom_{\f{g}} \big(T_{\tau+\nu}^{\tau+\nu+\mu}  T_{\tau}^{\tau+\nu}
                            \verma{x \cdot \tau},
                          T_{\tau}^{\tau+\nu+\mu} \verma{x \cdot \tau}\big)\\
    & \stackrel{\sim}{\longrightarrow}   
      \Hom_{\f{g}} \big(\verma{x \cdot (\tau+\nu+\mu)},
                           \verma{x \cdot (\tau+\nu+\mu)}\big)\\
    & \stackrel{\sim}{\longrightarrow} 
      k
  \end{align*}
  \end{center}

 Denote this map by $\Nat(\mu, \nu;x)(\tau)$.
 Here again the first isomorphism is clear by the Theorem of Bernstein-Gelfand,
   the second is due to the the maps $F_{x(\nu+\mu)}(x \cdot \tau)$ and 
   $F_{x \mu}(x\cdot(\tau+\mu))\circ 
                        T_{\tau +\nu}^{\tau +\nu+\mu} F_{x\nu}(x \cdot \tau)$.
 We then obtain
   \[ 
   \Delta(\mu, \nu; x)(\tau) = \Nat(\mu, \nu;x)(\tau) \circ 
                                        (\Nat(\mu, \nu;e)(\tau))^{-1} (1). 
   \] 
 It follows that this characterization of $\Delta$ is independent of the
  choice of the weight vectors $e_{\nu}, \tilde{e}_{\mu}$ and 
  $\bar{e}_{\mu +\nu}$, and also independent of the choice of the pre-image
  $\dot{x}$ of $x$.

     \section{\bf Bernstein's relative trace and a special case of
                                          $\boldsymbol \Delta$} \label{bern}

\subsection{The relative trace}

   Let us recall the definition of the 
      {\em relative trace} $\tr_{E}$ defined by Bernstein in \cite{be}.
   Let $\f{g}$-mod be the category of all $\f{g}$-modules and let $E$ be a
      finite dimensional $\f{g}$-module. 
   Denote by $F_E : \gmod \to \gmod$ the functor defined by
      $F_E (M) := E \otimes M$.
   The relative trace 
      $\tr_E : \End_{\gmod}(F_E) \to \End_{\gmod}(\Id)$
      is then a morphism from the endomorphisms of the functor $F_E$ to the 
      endomorphisms of the identity functor on $\gmod$, defined by
   \[
   \begin{array}{lccc}
         \tr_E^M : & \End_{\f{g}}(E \otimes M) & \to     & \End_{\f{g}} M\\
                   & a                         & \mapsto & \tr_E^M (a)
   \end{array}
   \]
   where
   \[\begin{CD}
     \tr^M_E (a): \;
                   M 
                   @>{\inc}>>
                   E^{*} \otimes E \otimes M 
                   @>{\id \otimes a}>>
                   E^{*} \otimes E \otimes M 
                   @>{\cont \otimes \id}>>
                   M 
    \end{CD} .\] 
  
   Here, $\inc$ is the map
                 $M 
                 \stackrel{\jott}{\hookrightarrow} 
                 \End_{\f{g}}(E) \otimes M
                 \stackrel{\ce}{\to}
                 E^* \otimes E \otimes M$
         with $\jott(m) := \id_E \otimes m$ and $\ce$ the canonical 
         isomorphism $\End_{\f{g}}(E) \cong E^* \otimes E$.
   The map
   $\cont : E^* \otimes E \to k$ denotes the evaluation map.

   Bernstein has calculated an explicit formula for the relative trace, 
     by considering it as a map from ${\cal P}(\f{h}^*)^{\cal W \cdot}$
     to itself in the following way:
   First, we identify $\End_{\gmod}(\Id) \cong \f{Z}$ with the center of 
     the enveloping algebra $\f{U}$, then we make use of the natural
     morphism $\f{Z} \to \End_{\gmod}(F_E)$ and composing this with the
     trace map we obtain $\tr_E : \f{Z} \to \f{Z}$.
   By means of the Harish-Chandra isomorphism 
     $\xi:\,\f{Z} \stackrel{\sim}{\longrightarrow}
                                      {\cal P}(\f{h}^*)^{\cal W \cdot}$ 
     (normalized by $z-\xi(z) \in \f{U}\f{n}$)
     we may then regard $\tr_E$ as an endomorphism of 
     ${\cal P}(\f{h}^*)^{\cal W \cdot}$.      

   Define now on ${\cal P}(\f{h}^*)$ a convolution $f \mapsto P(E) \ast f$
     by $(P(E) \ast f)(\lambda) := \sum_{\mu \in P(E)} f(\lambda +\mu)$,
     where the sum is taken over all weights $\mu \in P(E)$ with their 
     multiplicities.
   Set $\Lambda(\lambda) 
   := \prod_{\alpha \in R^+} \langle \lambda +\rho,\alpha^{\vee}\rangle$.
   Then we have

   \begin{theoohne}\cite{be}
    $\tr_E(f) = \Lambda^{-1} (P(E) \ast \Lambda f)$ for all 
      $f \in {\cal P}(\f{h}^*)^{\cal W \cdot}$.
   \end{theoohne}

   If we choose now for $M$ the Verma module $\verma{\lambda}$ we can
     associate to each endomorphism $f \in \End_{\f{g}}
     (E \otimes \verma{\lambda})$ an element
     $\tr_E^{M(\lambda)}(f)$ of $\End_{\f{g}}(\verma{\lambda}) \cong k$.       
   As endomorphism of $\verma{\lambda}$ this element operates on 
     $\verma{\lambda}$ by multiplication with the scalar 
     $(\tr_E^{M(\lambda)}(f))(\lambda)$ and we obtain by Bernstein's
     Theorem for all $\lambda \in \f{h}^*$ and for all $w \in \cal W$
    \begin{align*}
          \big(\tr_E^{M(w\cdot\lambda)}(f)\big)(\lambda) 
          &=
          \big(\Lambda^{-1} (P(E) \ast \Lambda f)\big)(\lambda)\\ 
          &=
          \left(\Lambda(\lambda)\right)^{-1}
                \sum_{\mu \in P(E)}\Lambda(\lambda +\mu) f(\lambda +\mu) .
    \end{align*} 

\subsection{The special case $\Delta(-\nu,\nu;w_0)$}

   Let now $E=E(\nu)$, $w_0 \in \cal W$ the longest element and
      $\pr_{\chi(\tau +\nu)} \in 
                        \End_{\f{g}}(E(\nu) \otimes \verma{w_0 \cdot \tau})$
      the projection on the central character $\chi(\tau +\nu)$.
   Then $\tr_{E(\nu)}^{M(w_0 \cdot \tau)} (\pr_{\chi(\tau +\nu)})$
      is an element in $\End_{\f{g}}(\verma{w_0 \cdot \tau}) \cong k$ 
      and we have

   \begin{theo} \label{Bernstein}
      Let $\nu \in P^+$ be a dominant weight and $\tau \in \f{h}^*$ generic.
      Then 
      \[ \Delta(-\nu, \nu; w_0)(\tau) = 
        \tr_{E(\nu)}^{M(w_0 \cdot \tau)} (\pr_{\chi(\tau +\nu)}) .\]
   \end{theo} 

   \begin{proof} Postponed to \ref{bew_theorem}. \end{proof}

   Regarding $\Delta(-\nu, \nu;w_0)$ as a polynomial function on $\f{h}^*$ 
     we obtain for this special case an explicit formula:

   \begin{kor}
      Let $\nu \in P^+$ be a dominant weight and $\tau \in \f{h}^*$ generic.
      Then 
      \[ \Delta(-\nu, \nu;w_0)(\tau) 
      =
      \prod_{\alpha \in R^+} \frac{\langle \tau +\nu +\rho, 
                                             \alpha^{\vee} \rangle}
                                  {\langle \tau +\rho, 
                                             \alpha^{\vee} \rangle} \; .\]
   \end{kor}

  \begin{proof}
    Under the morphism ${\cal P}(\f{h}^*)^{\cal W \cdot} \cong \f{Z} 
      \to \End_{\gmod}(F_E)$ the projection $\pr_{\chi(\tau +\nu)}$ 
      is the image of
      a polynomial function, which takes value $1$ at all weights 
      $\lambda \in {\cal W} \cdot (\tau +\nu)$ 
      and vanishes at all other weights $\tau +\mu$ with $\mu \in P(E)$.
    Call this polynomial function $\fpr$. 
    By Bernstein's formula for the relative trace we then obtain
     \[
       \big(\tr_{E(\nu)}^{M(w_0 \cdot \tau)}(\fpr)\big)(\tau) =
       (\Lambda (\tau))^{-1}
          \sum_{\mu \in P(E(\nu))} \Lambda (\tau +\mu) \fpr(\tau +\mu)
     \]
      and by definition of $\fpr$ the value of $\fpr(\tau +\mu)$ for  
      $\mu \in P(E)$ does not vanish if and only if there is a $w \in \cal W$
      such that $w \cdot(\tau +\mu) =(\tau +\nu)$.
    Since $\tau$ is generic, this is only possible for $w=e$ and hence 
      $\mu =\nu$.
    In this case we have $\fpr(\tau +\nu) = 1$ and since furthermore
      $\dim E(\nu)_{\nu} = 1$, we get
      $\sum_{\mu \in P(E(\nu))} \Lambda (\tau +\mu) \fpr(\tau +\mu)
       = \Lambda (\tau  + \nu)$.
    The claim now follows by Theorem \ref{Bernstein} and the equation
      \[\big(\tr_{E(\nu)}^{M(w_0 \cdot \tau)}(\fpr)\big)(\tau) 
        = \frac{\Lambda (\tau + \nu)}{\Lambda (\tau)}
        = \prod_{\alpha \in R^+} \frac{\langle \tau +\nu +\rho, 
                                             \alpha^{\vee} \rangle}
                                  {\langle \tau +\rho, 
                                             \alpha^{\vee} \rangle} \; .\]
     \end{proof}

\subsection{Proof of Theorem \ref{Bernstein}} \label{bew_theorem}

  First we make some more general preliminary remarks.

  \subsubsection{The adjunctions $(T_{\tau}^{\tau+\nu},T_{\tau+\nu}^{\tau})$ 
                 and $(T^{\tau}_{\tau+\nu},T^{\tau+\nu}_{\tau})$}
 
     Let $\cal A$ and $\cal B$ be categories and $F:\cal A \to \cal B$, 
        $G:\cal B \to \cal A$ two functors.
     Then an adjunction $(F,G)$ of $F$ and $G$ is a family of isomorphisms
        \[ (F,G)_{M,N} := (F,G) :
                      \Hom_{\cal B}(FM,N) 
                      \stackrel{\sim}{\longrightarrow} 
                      \Hom_{\cal A}(M,GN) , \]
        which is natural in $M$ and $N$ ($M \in \cal A$,\; $N \in \cal B$).

     For example, for $E$ a finite dimensional $\f{g}$-module we obtain
       an adjunction $(F_E,F_{E^*})$ of $F_E$: $\gmod \to \gmod$
       and $F_{E^*}$: $\gmod \to \gmod$ as the composition
       \[\Hom_{\f{g}}(E \otimes M,N)   
         \to
         \Hom_{\f{g}}(E^* \otimes E \otimes M,E^* \otimes N)
         \to
         \Hom_{\f{g}}(M,E^* \otimes N) .
        \]   
     Here, the first map is given by $f \mapsto \id_{E^*} \otimes f$,
       the second by $g \mapsto g \circ \inc$. 
     Interchanging $E$ and $E^*$ we obtain in the same way an adjunction
       $(F_{E^*}, F_E)$. 
     Its inverse is the composition
       \[\Hom_{\f{g}}(M,E \otimes N)    
         \to 
         \Hom_{\f{g}}(E^* \otimes M,E^* \otimes E \otimes N) 
         \to
         \Hom_{\f{g}}(E^* \otimes M,N) , 
       \]
       where the first map is again given by $f \mapsto \id_{E^*} \otimes f$ 
       and the second by $g \mapsto (\cont \otimes \id_N) \circ g$.
     Let now $\nu \in P$ be an integral weight.
     Each identification
       $\varphi : E(\nu)^* \stackrel{\sim}{\longrightarrow} E(-\nu)$ 
       then defines adjunctions $(F_{E(\nu)}, F_{E(-\nu)})$ and 
       $(F_{E(-\nu)}, F_{E(\nu)})$.
     More precisely, we have
       \[
       \begin{array}{lccc}
         (F_{E(\nu)}, F_{E(-\nu)}) : \, &
               \Hom_{\f{g}}(E(\nu) \otimes M,N) &
               \stackrel{\sim}{\longrightarrow} &
               \Hom_{\f{g}}(M,E(-\nu) \otimes N) \\
               & f & \mapsto & (\varphi \otimes f) \circ \inc
       \end{array}
       \]
       and
       \[
       \begin{array}{lccc}
         (F_{E(-\nu)}, F_{E(\nu)})^{-1} :\,&
            \Hom_{\f{g}}(M,E(\nu) \otimes N) &
            \stackrel{\sim}{\longrightarrow} &
            \Hom_{\f{g}}(E(-\nu) \otimes M,N) \\
            & g& \mapsto & (\cont \otimes \id_N) \circ (\varphi^{-1} \otimes g)
        \end{array}
        \]
     Let now $\inc_{\chi} :\cal M^{\infty}(\chi) \hookrightarrow \cal M$
       denote the embedding functor.
     We then have in a natural way adjunctions $(\inc_{\chi},\pr_{\chi})$ 
       and $(\pr_{\chi},\inc_{\chi})$.
     Since for the translation functor 
        $T_{\tau}^{\tau+\nu} =
            \pr_{\chi(\tau+\nu)} \circ F_{E(\nu)} \circ \inc_{\chi(\tau)}$,
        we thus obtain also adjunctions 
        $(T_{\tau}^{\tau+\nu},T^{\tau}_{\tau+\nu})$ and 
        $(T^{\tau}_{\tau+\nu},T_{\tau}^{\tau+\nu})$.

\subsubsection{The natural transformations $\adj^1$ and $\adj^2$} 
  \label{ber}

     Let $M \in \cal M^{\infty}(\chi(\tau))$ and consider 
        $\Id \in \Hom_{\f{g}}(T_{\tau}^{\tau+\nu} M,T_{\tau}^{\tau+\nu} M)$.
     By means of the two adjunctions 
        $(T_{\tau}^{\tau+\nu},T_{\tau+\nu}^{\tau})$ and
        $(T^{\tau}_{\tau+\nu},T^{\tau+\nu}_{\tau})^{-1}$ we get two
        maps as the images of $\Id$:
        \[\adj^1_M \in 
             \Hom_{\f{g}}(M,T_{\tau+\nu}^{\tau}  T_{\tau}^{\tau+\nu} M)
        \text{ and }
        \adj^2_M \in
             \Hom_{\f{g}}(T_{\tau+\nu}^{\tau}  T_{\tau}^{\tau+\nu} M,M)\]
        and one checks that we obtain in this way natural transformations
        $\adj^1 \in 
             \Homtauinf(\Id,T_{\tau+\nu}^{\tau}\circ T_{\tau}^{\tau+\nu})$
        and
        $\adj^2 \in
             \Homtauinf(T_{\tau+\nu}^{\tau} \circ T_{\tau}^{\tau+\nu},\Id)$.
     By composing these two maps, we get a canonical endomorphism of $M$:
        \[
        \begin{CD}
                    M     @>{\adj^1_M}>>  
                    T_{\tau+\nu}^{\tau}  T_{\tau}^{\tau+\nu} M
                    @>{\adj^2_M}>>
                    M .
        \end{CD}
        \]
     We want to describe this endomorphism in more detail.
     By definition, the adjunction $(T_{\tau}^{\tau+\nu},T^{\tau}_{\tau+\nu})$
       is just the composition of adjunctions
       $(\inc_{\chi(\tau)},\pr_{\chi(\tau)}) \circ
       (F_{E(\nu)}, F_{E(-\nu)}) \circ 
       (\pr_{\chi(\tau+\nu)},\inc_{\chi(\tau+\nu)})$
       and one checks easily that the image of 
       $\Id \in \Hom_{\f{g}}(T_{\tau}^{\tau+\nu}M,T_{\tau}^{\tau+\nu}M)$ under
       $(\pr_{\chi(\tau+\nu)},\inc_{\chi(\tau+\nu)})$ 
       is precisely the projection $\pr_{\chi(\tau+\nu)}$.
     Thus, the image of $\Id$ under the adjunction  
       $(T_{\tau}^{\tau+\nu},T^{\tau}_{\tau+\nu})$ can be described by the 
       composition
        \[\adj^1_M :
        \begin{CD}
             M
             @>{\inc}>>
             E^* \otimes E \otimes M
             @>{\varphi \otimes \pr}>>
             T_{\tau+\nu}^{\tau}  T_{\tau}^{\tau+\nu} M.
           \end{CD}
         \]
     Here, we put $E = E(\nu)$ and $\pr =\pr_{\chi(\tau+\nu)}$. 
     Analogously, we obtain for $\adj^2_M$ the composition
         \[\adj^2_M :
         \begin{CD}
              T_{\tau+\nu}^{\tau} T_{\tau}^{\tau+\nu} M  
              @>{\varphi^{-1} \otimes \pr}>>
              E^* \otimes E \otimes M
              @>{\cont \otimes \id_M}>>
              M 
         \end{CD} 
         \]
       and taken together we have the following commutative diagram:
         \[
         \begin{CD}
              M
              @>{\adj^1_M}>>
              T_{\tau+\nu}^{\tau}  T_{\tau}^{\tau+\nu} M
              @>{\adj^2_M}>>
              M\\
              @V{\inc}VV
              @A{\varphi \otimes \id_{E \otimes M}}AA
              @A{\cont \otimes \id_M}AA\\
              E^* \otimes E \otimes M
              @>{\id_{E^*} \otimes \pr}>>
              E^* \otimes E \otimes M
              @>{\id_{E^*} \otimes \pr}>>
              E^* \otimes E \otimes M
          \end{CD}
          \] 
     We thus obtain for $\adj^2_M \circ \adj^1_M$
         \[
         \begin{CD}
                      M
                      @>{\inc}>>
                      E^* \otimes E \otimes M
                      @>{\id_{E^*} \otimes \pr_{\chi(\tau+\nu)}}>>  
                      E^* \otimes E \otimes M 
                      @>{\cont \otimes \id_M}>>
                      M.
          \end{CD}
          \]
    Comparing this with the relative trace $\tr_E$ for $E = E(\nu)$,
      it follows immediately for all $M \in \cal M^{\infty}(\chi(\tau))$ that
      \[\adj^2_M \circ \adj^1_M = \tr_{E(\nu)}^M(\pr_{\chi(\tau +\nu)}) \in
           \End_{\f{g}}(M)\]
      or, regarded as natural transformation of the identity functor
      $\Id : \cal M^{\infty}(\chi(\tau))$ $\to \cal M^{\infty}(\chi(\tau))$
      to itself
      \[\adj^2 \circ \adj^1 = \tr_{E(\nu)}(\pr_{\chi(\tau +\nu)}) .\]

    Let now $\iota_{\chi} : {\cal M}(\chi) \hookrightarrow \cal M$ denote the 
      embedding functor, let $F, G : {\cal M}^{\infty}(\chi) \to \cal M$ 
      be functors and denote by $F(\chi)$, resp. $G(\chi)$ its restrictions 
      to the subcategory ${\cal M}(\chi)$.
    Each natural transformation $n$ from $F$ to $G$ can be regarded as
      natural transformation from $F(\chi)$ to $G(\chi)$ by first applying
      $\iota_{\chi}$ and then $n$. 
    In particular, we obtain in this way the two natural transformations
      $\adj^1 \in 
      \Homtau \big(\Id(\chi), T_{\tau+\nu}^{\tau}\circ 
                                          T_{\tau}^{\tau+\nu}(\chi)\big)$
      and similarly
      $\adj^2 \in 
      \Homtau \big(T_{\tau+\nu}^{\tau} \circ T_{\tau}^{\tau+\nu}(\chi),
                                          \Id(\chi)\big)$. 

 \subsubsection{We have $\Delta(-\nu, \nu;w_0)(\tau) = 
                 \adj^2_{M(w_0 \cdot \tau)}\circ \adj^1_{M(w_0 \cdot \tau)}$}

    By choosing for $M$ the Verma module $\verma{\tau}$, we can assign to
      each $\tau \in \f{h}^*$ a canonical element  
      $\adj^2_{M(\tau)} \circ \adj^1_{M(\tau)} \in \End(\verma{\tau}) \cong k$.
    We will see in the following, that for dominant $\nu$ this is precisely 
      the triangle function $\Delta(-\nu, \nu;w_0)(\tau)$.
    For $x \in \cal W$ and $\tau$ generic we call $\F_x(x \cdot\tau)$ 
      the isomorphism
      $T_{\tau+\nu}^{\tau}  T_{\tau}^{\tau+\nu} \verma{x \cdot \tau}
      \stackrel{\sim}{\longrightarrow} 
      \verma{x \cdot \tau}$ 
      given by  
      $(\F_x(x \cdot \tau))^{-1}= 
       T_{\tau +\nu}^{\tau} F_{x \nu}(x \cdot \tau)\circ 
                                           F_{- x \nu}(x \cdot (\tau +\nu))$ 
      (see \ref{uncan}).
    First we show that for dominant $\nu$ and generic $\tau$ the following
      diagram commutes:
      \[\xymatrix{
            {T_{\tau+\nu}^{\tau} T_{\tau}^{\tau+\nu} \verma{\tau}}
            \ar[rr]^-{\adj^2_{M(\tau)}}
            \ar[d]_{\F_e(\tau)} &&
            {\verma{\tau}}\\
            {\verma{\tau}} \ar[urr]_{\id}
          }
      \]
    Let again denote $e_{\nu} \in E(\nu)_{\nu}$ the fixed extremal weight
      vector for a dominant integral weight $\nu$.
    Since $w_0$ is the longest element in the Weyl group the dominance of
      $\nu$ implies that $e_{-\nu} := \dot{w}_0 e_{\nu}$ is a weight vector
      of weight $-\nu$.
    Here, $\dot{w}_0 \in G$ is a representative of $w_0 \in N_G(T)/T$. 

    Define now a pairing $E(-\nu) \times E(\nu) \to k$ by 
      $\langle e_{-\nu},e_{\nu} \rangle := 1$ and obtain thus an identification
      $\varphi : E(\nu)^* \stackrel{\sim}{\longrightarrow} E(-\nu)$.
    Let $v_{\tau} \in \verma{\tau}$ be the canonical generator.
    To see that the above diagram commutes, it suffices to show that the 
      pre-image of $v_{\tau}$ in 
      $T_{\tau+\nu}^{\tau} T_{\tau}^{\tau+\nu} \verma{\tau}$
      is mapped again to $v_{\tau}$ when applying $\adj^2_{M(\tau)}$.
    We have
      \begin{align*}
           (\F_e(\tau))^{-1}(v_{\tau}) 
           &=
           \pr_{\chi(\tau)}\big( e_{-\nu} \otimes (\pr_{\chi(\tau+\nu)}
                                        (e_{\nu} \otimes v_{\tau}))\big)\\ 
           &=
           \pr_{\chi(\tau)}( e_{-\nu} \otimes e_{\nu} \otimes v_{\tau}) \\
           &=
           e_{-\nu} \otimes e_{\nu} \otimes v_{\tau}
              - {\bigoplus}_{\chi \neq \chi(\tau)}  
                \pr_{\chi} (e_{-\nu} \otimes e_{\nu} \otimes v_{\tau}).
         \end{align*}

    The first equation holds by definition of $\F_e(\tau)$, the second
      follows since for dominant $\nu$ we have
      $\pr_{\chi(\tau +\nu)}(e_{\nu} \otimes v_{\tau}) 
                                           = e_{\nu} \otimes v_{\tau}$
      (see the example on page \pageref{beispiel})
      and the third equation is just direct sum decomposition.
    Since $\adj^2_M$ was the composition
      \[ 
      \begin{CD}
              T_{\tau+\nu}^{\tau} T_{\tau}^{\tau+\nu} M  
              @>{\varphi^{-1} \otimes \pr}>>
              E^* \otimes E \otimes M
              @>{\cont \otimes \id_M}>>
              M 
       \end{CD} 
       \]
      we know in particular that $\adj^2_{M(\tau)}$ is a $\f{g}$-module 
      homomorphism to $\verma{\tau}$. 
    Then it is clear, that $\adj^2_{M(\tau)}$ maps 
      ${\bigoplus}_{\chi \neq \chi(\tau)}  
                \pr_{\chi} (e_{-\nu} \otimes e_{\nu} \otimes v_{\tau})$ to 
      zero since this element has wrong central character.
    For $e_{-\nu} \otimes e_{\nu} \otimes v_{\tau}$ we use the fact that
      $\langle e_{-\nu},e_{\nu} \rangle = 1$ and obtain as image under
      $\adj^2_{M(\tau)}$ precisely $v_{\tau}$.
    Taken together, we get $\adj^2_{M(\tau)} = \id \circ \F_e(\tau)$, i.e. the
      above diagram commutes.

    This now means that under the map $\Nat(-\nu,\nu;e)(\tau)$ the image of 
      $\adj^2 \in \Homtau \big(T_{\tau+\nu}^{\tau}\circ 
                       T^{\tau+\nu}_{\tau}(\chi), T_{\tau}^{\tau}(\chi)\big)$
      is just the identity $\id \in \End(\verma{\tau})$.
    By means of the uncanonical definition of the triangle function
      we deduce
      \begin{align*}
            \Delta(-\nu, \nu; w_0)(\tau)
            &=
            \Nat(-\nu,\nu;w_0)(\tau)\circ \big(\Nat(-\nu,\nu;e)(\tau)\big)^{-1}
                                                           (\id_{M(\tau)})\\
            &= \big(\Nat(-\nu,\nu;w_0)(\tau)\big)(\adj^2)
      \end{align*}
      and the Theorem of Bernstein-Gelfand implies that the right hand side
      of the following diagram commutes:\\

      \begin{minipage}{1cm}
          \bf ($\ast$)
      \end{minipage}
      \begin{minipage}{7cm}
      $ 
         \xymatrix{
         {\verma{w_0 \cdot \tau}}
         \ar[rr]^-{\adj^1_{M(w_0 \cdot \tau)}}
         \ar[drr]_{\id} &&
         {T_{\tau+\nu}^{\tau} T_{\tau}^{\tau+\nu}\verma{w_0 \cdot \tau}}
         \ar[rr]^-{\adj^2_{M(w_0 \cdot \tau)}}
         \ar[d]^{\F_{w_0}(w_0 \cdot \tau)} &&
         {\verma{w_0 \cdot \tau}}\\ 
         &&{\verma{w_0 \cdot \tau}} 
         \ar[urr]_-{\;\;\;\;\;\;\;\;\;\Delta(-\nu, \nu;w_0)(\tau)}
          } 
        $
       \end{minipage}
 
     The commutativity of the left hand side can be shown in an analogous way.
     We thus obtain for dominant $\nu$
       \[\adj^2_{M(w_0 \cdot \tau)} \circ \adj^1_{M(w_0 \cdot \tau)}
         = \Delta(-\nu, \nu;w_0)(\tau)\]
       and together with the equation 
       $\adj^2_{M(w_0 \cdot \tau)} \circ \adj^1_{M(w_0 \cdot \tau)} 
          = \tr_{E(\nu)}^{M(w_0 \cdot \tau)}(\pr_{\chi(\tau +\nu)})$
       from section \ref{ber} we deduce Theorem \ref{Bernstein}.
   \qed

\subsection{The case $\verma{\tau} \longrightarrow  
         T_{\tau+\nu}^{\tau}  T_{\tau}^{\tau+\nu} \verma{\tau} \longrightarrow
         \verma{\tau}$}

  \begin{korbew}
      Let $\nu \in P^+$ be dominant and $\tau \in \f{h}^*$ generic.
      Then the composition    
      \[\begin{CD}
         \verma{\tau}  @>{\adj^1_{M(\tau)}}>>  
         T_{\tau+\nu}^{\tau}  T_{\tau}^{\tau+\nu} \verma{\tau}
                       @>{\adj^2_{M(\tau)}}>>
         \verma{\tau}
        \end{CD}
        \]
      is just multiplication with 
      $s(\tau) := \prod_{\alpha \in R^+} 
                  \frac{\langle \tau +\nu +\rho, \alpha^{\vee} \rangle}
                       {\langle \tau +\rho, \alpha^{\vee} \rangle}$, 
      where $s$ is considered as an element in $\Quot(S(\f{h}))$.
   \end{korbew}

   \begin{proof}
      The commutativity of the diagram $(\ast)$ in the proof of Theorem
         \ref{Bernstein} is equivalent to the commutativity of
         \[  
           \xymatrix{
           {\verma{\tau}}
            \ar[rr]^-{\adj^1_{M(\tau)}}
            \ar[drr]_-{\Delta(-\nu, \nu;w_0)(\tau)\;\;\;\;\;\;} &&
            {T_{\tau+\nu}^{\tau} T_{\tau}^{\tau+\nu}\verma{\tau}}
            \ar[rr]^-{\adj^2_{M(\tau)}}
            \ar[d]^{\F_{e}(\tau)} &&
            {\verma{\tau}}\\ 
            &&{\verma{\tau}} 
            \ar[urr]_{\id}
          } 
         \]
      Hence $(\adj^2 \circ \adj^1)_{M(\tau)} 
               = \Delta(-\nu, \nu;w_0)(\tau)
               = \prod_{\alpha \in R^+} 
                  \frac{\langle \tau +\nu +\rho, \alpha^{\vee} \rangle}
                       {\langle \tau +\rho, \alpha^{\vee} \rangle}
               = s(\tau)$.
     \end{proof}

     \section{\bf Calculation of $\boldsymbol \Delta$}\label{formel}

   \begin{theo} \label{dreiecksfunktion}
      For $\lambda \in \f{h}^*$ set $\bar{\alpha}(\lambda) := 1$ if
      $\langle\lambda,\alpha^{\vee}\rangle < 0$ and 
      $\bar{\alpha}(\lambda) := 0$ if 
      $\langle\lambda,\alpha^{\vee}\rangle \ge 0$.
      Let $\nu, \mu \in P$ be integral weights and $x \in \cal W$. 
      Then there exists a constant $c \in k^{\times}$, independent of 
      $\tau, \nu, \mu$ and $x$, such that
      \[\Delta(\mu, \nu; x)(\tau-\rho)
      = c \!\!\!\!\!\!\prod_{\genfrac{}{}{0pt}{2}
                {\alpha \in R^+}{\text{with } x\alpha \in R^-}}\!\!\!\!
        \frac
      {\langle \tau,\alpha^{\vee}\rangle^{\bar{\alpha}(\nu)}
       \langle\tau+\nu,\alpha^{\vee}\rangle^{\bar{\alpha}(\mu)}
       \langle\tau+\nu+\mu,\alpha^{\vee}\rangle^{\bar{\alpha}(\nu+\mu)}}
      {\langle \tau,\alpha^{\vee}\rangle^{\bar{\alpha}(\nu+\mu)}
       \langle \tau+\nu,\alpha^{\vee}\rangle^{\bar{\alpha}(\nu)}
       \langle \tau+\nu+\mu,\alpha^{\vee}\rangle^{\bar{\alpha}(\mu)}}.\]
   \end{theo}

   \begin{proof} Postponed. \end{proof}

   For an integral weight $\nu \in P$ define the map 
       $\delta_{\nu} \in {\cal P}(\f{h}^*)$ as in Section \ref{pole} by 
       \[\delta_{\nu}(\tau) :=
       \prod_{\alpha\in R^+}
       \;\;
       \prod_{0 \le n_{\alpha} < -\langle \nu, \alpha^{\vee}\rangle}
       (\langle \tau +\rho,\alpha^{\vee} \rangle -n_{\alpha}) \;.
       \]

   \begin{lemma} \label{langeformel}
      Let $\pi(\tau)$ be the product on the right hand side of the equation
        in Theorem \ref{dreiecksfunktion}. 
      Then
     \[\pi(\tau+\rho) = \pm \frac{\delta_{\nu}(\tau)\delta_{\mu}(\tau +\nu)
                                       \delta_{x(\nu+\mu)}(x \cdot\tau)}
               {\delta_{x\nu}(x\cdot\tau)\delta_{x\mu}(x\cdot(\tau +\nu))
                                       \delta_{\nu+\mu}(\tau)} .\] 
   \end{lemma}
   
   \begin{proof} 
      Let $x \in \cal W$ be fixed and let $\alpha \in R^+$ be a positive root.
      We then have in $\delta_{\nu}(\tau)$ the factor
        $\prod_{0\le n < -\langle\nu,\alpha^{\vee}\rangle}
            (\langle\tau +\rho,\alpha^{\vee}\rangle -n)$.
      Suppose now that $x\alpha$ is also a positive root.
      Then we obtain in $\delta_{x\nu}(x \cdot\tau)$ the factor 
      \begin{align*}
        \prod_{0\le n < -\langle x\nu, x\alpha^{\vee}\rangle}\!\!\!\!
         (\langle x \cdot \tau +\rho, x\alpha^{\vee}\rangle -n)
        &=
        \prod_{0\le n < -\langle \nu, \alpha^{\vee}\rangle}
         (\langle x (\tau +\rho), x\alpha^{\vee}\rangle -n)\\
        &=
        \prod_{0\le n < -\langle \nu, \alpha^{\vee}\rangle}
         (\langle \tau +\rho, \alpha^{\vee}\rangle -n)
       \end{align*}

       and hence in the quotient 
       $\delta_{\nu}(\tau)/\delta_{x\nu}(x \cdot\tau)$ all the products over
       $\alpha \in R^+$ with $x \alpha \in R^+$ cancel out.
     Let now $\alpha \in R^+$ such that $x \alpha \notin R^+$, i.e.
       $-x \alpha$ is a positive root. 
     In this case we have in $\delta_{x\nu}(x \cdot\tau)$ the factor
       \begin{align*}
        \prod_{0\le n < -\langle x\nu, -x\alpha^{\vee}\rangle}\!\!\!\!
         (\langle x \cdot \tau +\rho, -x\alpha^{\vee}\rangle -n)
        &=
        \prod_{0\le n < \langle \nu, \alpha^{\vee}\rangle}
         (-\langle \tau +\rho, \alpha^{\vee}\rangle -n)\\
        &=
        (-1)^{\langle \nu, \alpha^{\vee}\rangle}\!\!\!\!
        \prod_{0\le n < \langle \nu, \alpha^{\vee}\rangle}
         \!\!\!\!(\langle \tau +\rho, \alpha^{\vee}\rangle +n)
       \end{align*}
        
     and taken together we get
      \[\frac{\delta_{\nu}(\tau)}{\delta_{x\nu}(x \cdot\tau)}
        =
        \prod_{\genfrac{}{}{0pt}{2}
             {\alpha \in R^+}{\text{with } x\alpha \in R^-}}\!\! \!\!
         (-1)^{\langle \nu, \alpha^{\vee}\rangle}
         \frac
         {\prod_{0\le n_{\alpha} < -\langle \nu, \alpha^{\vee}\rangle}
          (\langle \tau +\rho, \alpha^{\vee}\rangle -n_{\alpha})} 
         {\prod_{0\le n_{\alpha} < \langle \nu, \alpha^{\vee}\rangle}
           (\langle \tau +\rho, \alpha^{\vee}\rangle +n_{\alpha})} .\]

      Considering then the different cases according to whether
          $\langle \nu, \alpha^{\vee}\rangle , 
           \langle \mu, \alpha^{\vee}\rangle$ and
          $\langle \nu + \mu, \alpha^{\vee}\rangle$ are positive or
          negative, we obtain the closed formula
          \[\frac{\delta_{\nu}(\tau)\delta_{\mu}(\tau +\nu)
                                       \delta_{x(\nu+\mu)}(x \cdot\tau)}
           {\delta_{x\nu}(x\cdot\tau)\delta_{x\mu}(x\cdot(\tau +\nu))
                                       \delta_{\nu+\mu}(\tau)}
          = \varepsilon \; \pi(\tau+\rho) ,\]  
          where $\varepsilon := \prod\limits_{\genfrac{}{}{0pt}{2}
             {\alpha \in R^+}{\text{mit } x\alpha \in R^-}}\!\! \!\!
     (-1)^{(\bar{\alpha}(\nu) \langle \nu,\alpha^{\vee}\rangle
          +\bar{\alpha}(\mu) \langle \mu,\alpha^{\vee}\rangle
          -\bar{\alpha}(\nu+\mu)\langle \mu+\nu,\alpha^{\vee}\rangle)}
     =\pm 1. $ \qed
   \renewcommand{\qed}{}
   \end{proof}

   Now, let $\tau \in \f{h}^*$ be a generic weight and recall 
   (see \ref{uncan}) the map
     $\Nat_x(\tau) := \Nat(\mu,\nu;x)(\tau)$ :\\

  \hspace*{0.5em}
  $\Homtau \big(T_{\tau +\nu}^{\tau +\nu +\mu} \circ T_{\tau}^{\tau +\nu},
                               T_{\tau}^{\tau +\nu +\mu}\big)$ \\[-2.5em]

  \begin{center}
  \begin{align*}
    & \stackrel{\sim}{\longrightarrow} 
      \Hom_{\f{g}} \big(T_{\tau+\nu}^{\tau+\nu+\mu}  T_{\tau}^{\tau+\nu}
                            \verma{x \cdot \tau},
                          T_{\tau}^{\tau+\nu+\mu} \verma{x \cdot \tau}\big)\\
    & \stackrel{\sim}{\longrightarrow}   
      \Hom_{\f{g}} \big(\verma{x \cdot (\tau+\nu+\mu)},
                           \verma{x \cdot (\tau+\nu+\mu)}\big)\\
    & \stackrel{\sim}{\longrightarrow} 
      k
  \end{align*}
  \end{center}

  Denote the pre-image of $1 \in k$ under this isomorphism by the natural
    transformation
    $ G^x(\tau) := (\Nat_x(\tau))^{-1}(1) \in 
    \Homtau (T_{\tau +\nu}^{\tau +\nu +\mu} \circ T_{\tau}^{\tau +\nu},
                               T_{\tau}^{\tau +\nu +\mu})$.

  \begin{lemma} \label{gleichung}
      For generic $\tau \in \f{h}^*$ we have
         $G^e(\tau)= \Delta(\mu,\nu;x)(\tau) G^x(\tau)$.
  \end{lemma}

  \begin{proof}
     By definition,
     $\Delta(\mu,\nu;x)(\tau)=\Nat_x(\tau)\circ(\Nat_e(\tau))^{-1} (1)$
     and therefore 
     \begin{align*}
       G^e(\tau)& =(\Nat_x(\tau))^{-1}\circ
                                    \Nat_x(\tau)\circ(\Nat_e(\tau))^{-1}(1)\\
                & =(\Nat_x(\tau))^{-1}(\Delta(\mu,\nu;x)(\tau))\\
                & =\Delta(\mu,\nu;x)(\tau)(\Nat_x(\tau))^{-1}(1)\\
                & =\Delta(\mu,\nu;x)G^x(\tau).      
     \end{align*}
  \end{proof} 

   Let now integral weights $\nu_1, \ldots, \nu_n$ be given such that
     $\sum_{i=1}^n \nu_i =0$.
   We set the translation functor
     $T(\nu_i):= T_{\tau+\nu_1+\cdots +\nu_{i-1}}^{\tau+\nu_1+\cdots +\nu_i}$ 
     when there is no ambiguity of the respective categories.
   For $\lambda \in P$ an integral weight and $\tau$ generic there are 
     isomorphisms (see \ref{uncan})
     $F_{x\lambda}(x \cdot \tau) : 
       \verma{x \cdot(\tau +\lambda)} \stackrel{\sim}{\longrightarrow} 
                               T_{\tau}^{\tau +\lambda} \verma{x \cdot \tau}$, 
     such that $(F_{x\lambda}(x \cdot \tau))(v_{x\cdot\tau}) =
        \pr_{\chi(\tau+\lambda)}(\dot{x} e_{\lambda} \otimes v_{x\cdot\tau})$.
   Here, $e_{\lambda} \in E(\lambda)_{\lambda}$ is a fixed chosen extremal 
     weight vector.
   Since $\sum_{i=1}^n \nu_i =0$, we can compose these isomorphisms to obtain
     an isomorphism
     $\verma{x\cdot \tau} \cong T(\nu_n)\cdots T(\nu_1)\verma{x \cdot \tau}$
     and thus also 
  \begin{align*}
    \Hom_{\f{g}} \big(\verma{x \cdot \tau},\verma{x \cdot \tau}\big)
    & \stackrel{\sim}{\longrightarrow} 
      \Hom_{\f{g}} \big(\verma{x \cdot \tau},
               T(\nu_n)\cdots T(\nu_1)\verma{x \cdot \tau}\big)\\
     & \hookrightarrow
       \Hom_k\big(\f{U}(\f{n}^-), 
           E(\nu_n)\otimes \cdots \otimes E(\nu_1)\otimes \f{U}(\f{n}^-)\big),
  \end{align*}
 
     where we have identified the Verma module $\verma{x\cdot \tau}$ with
     $\f{U}(\f{n}^-)$.
  Call now the image of the identity on $\verma{x\cdot \tau}$ under this map 
     $h^x(\nu_1, \ldots, \nu_n)(\tau)$ and let $\cal U \subset \f{h}^*$ 
     be the set of all generic weights.
  We then obtain for all $x \in \cal W$ a function
    \[
    \begin{array}{lccc}
      h^x(\nu_1,\ldots,\nu_n):\, 
       &\cal U 
       &\to 
       &\Hom_k \big(\f{U}(\f{n}^-),E(\nu_n)\otimes \cdots \otimes E(\nu_1)
                                             \otimes\f{U}(\f{n}^-)\big)\\
       & \tau    
       & \mapsto         
       & h^x(\nu_1,\ldots,\nu_n)(\tau),
      \end{array}
      \] 
     such that $h^x(\nu_1,\ldots,\nu_n)(\tau)$ maps the element
     $v_{x\cdot \tau} \in \verma{x\cdot \tau} \cong \f{U}(\f{n}^-)$ to the
     element $\pr_{\chi(\tau+\nu_1+\cdots +\nu_n)}(\dot{x}e_{\nu_n}\otimes 
     (\cdots \otimes
      \pr_{\chi(\tau+\nu_1)}(\dot{x}e_{\nu_1}\otimes v_{x\cdot \tau}))\cdots)$
     for fixed vectors $e_{\nu_i} \in E(\nu_i)_{\nu_i}$.

  \begin{lemma}\label{faktor}
     Set $d^x(\nu_1,\ldots,\nu_n)(\tau) :=
         \delta_{x\nu_1}(x\cdot\tau)\delta_{x\nu_2}(x\cdot(\tau+\nu_1))
          \cdots \delta_{x\nu_n}(x\cdot(\tau+\nu_1+\cdots +\nu_{n-1})).$
     Then the map
         $d^x(\nu_1,\ldots,\nu_n) h^x(\nu_1,\ldots,\nu_n)$ is algebraic on
         $\cal U$ and there exists an algebraic extension on $\f{h}^*$ 
         whose set of zeros has codimension $\ge 2$.
  \end{lemma}

  \begin{bemerkung}
    Here, we call a map $a:\; {\cal U} \to W$ to a vector space $W$
       {\em algebraic}, if it is a morphism of varieties.
    Of course, this is defined only if $\dim W < \infty$.
    In our case however, the image of $h^x(\nu_1, \ldots, \nu_n)$ is
       always contained in a finite dimensional subspace of
       $\Hom_k \big(\f{U}(\f{n}^-),E(\nu) \otimes \cdots \otimes E(\nu_1) 
                                                  \otimes \f{U}(\f{n}^-)\big)$
       and we may thus regard $h^x(\nu_1, \ldots, \nu_n)$ as a map
       between varieties.
  \end{bemerkung}

  \begin{proof} 
    Let $\nu \in P$ be an integral weight and recall the maps 
      \[
      \begin{array}{lccc}
        f_\nu:\, &\f{h}^* 
            & \longrightarrow 
            & E(\nu) \otimes \verma{\tau} \cong E(\nu) \otimes \f{U}(\f{n}^-)\\
            & \tau   
            & \mapsto         
            & \pr_{\chi(\tau+\nu)}(e_{\nu} \otimes v_{\tau})
       \end{array}
       \]       

     Let now $x \in \cal W$ be fixed.
     For generic $\tau$ define then the map $a^x_{\nu}(\tau)$ by
      \[
      \begin{array}{lccc}
        a^x_{\nu}(\tau):\, &\verma{x\cdot \tau} \cong \f{U}(\f{n}^-)
            & \longrightarrow 
            & E(\nu) \otimes \verma{x\cdot\tau} 
                                   \cong E(\nu)\otimes \f{U}(\f{n}^-)\\
            & v_{x\cdot\tau}   
            & \mapsto         
            & \delta_{x\nu}(x\cdot \tau) f_{x\nu}(x\cdot \tau)
       \end{array}
       \]       
  
       where $f_{x\nu}(x\cdot \tau) = 
             \pr_{\chi(\tau+\nu)}(\dot{x}e_{\nu} \otimes v_{x\cdot \tau})$.
     Since $f_{\nu}$ is algebraic on $\cal U$ (see Theorem \ref{theoeins}),
       we obtain in this way also an algebraic map
       \[
       \begin{array}{lccc}
        a^x_{\nu} :\, &\cal U 
         & \longrightarrow 
         & \Hom_k \big(\f{U}(\f{n}^-), E(\nu)\otimes\f{U}(\f{n}^-)\big)\\
         & \tau   
         & \mapsto         
         & a^x_{\nu}(\tau)
       \end{array}
       \] 
    
    Here again, the image of $a^x_{\nu}$ is contained in a finite dimensional
       subspace of 
       $\Hom_k \big(\f{U}(\f{n}^-),E(\nu) \otimes \f{U}(\f{n}^-)\big)$
       and we regard $a^x_{\nu}$ in this way as a map between varieties.

    According to Theorem \ref{theoeins} there is an algebraic extension of
       $\delta_{\nu}f_{\nu}$ on $\f{h}^*$, which vanishes only on a set of
       codimension $\ge 2$.
    Therefore, also $a^x_{\nu}$ has such an algebraic extension and we call
       it again $a^x_{\nu}$.
    For generic $\tau$ we then have 
       $(a^x_{\nu}(\tau))(v_{x\cdot \tau}) = \delta_{x\nu}(x\cdot \tau)
            \pr_{\chi(\tau+\nu)}(\dot{x}e_{\nu} \otimes v_{x\cdot \tau})$.
    Note that for all $\tau \in \f{h}^*$ the vector 
      $(a^x_{\nu}(\tau))(v_{x\cdot \tau}) \in E(\nu)\otimes\verma{x\cdot\tau}$
      generates a Verma module with highest weight $x \cdot(\tau+\nu)$.
    The image of $a^x_{\nu}(\tau)$ is thus always contained in a Verma module
      $\verma{x \cdot(\tau+\nu)} \subset E(\nu)\otimes \verma{x\cdot\tau}$
      and we can identify the element $(a^x_{\nu}(\tau))(v_{x\cdot \tau})$ 
      with the canonical generator
      $v_{x\cdot(\tau+\nu)} \in \verma{x \cdot(\tau+\nu)}$.
   Using the isomorphism $\verma{x \cdot(\tau+\nu)}\cong \f{U}(\f{n}^-)$, 
     we may then apply the map $a^x_{\nu^{\prime}}(\tau+\nu)$
     for another weight $\nu^{\prime} \in P$.
   In particular, we can concatenate the maps
     $a^x_{\nu_n}(\tau+\nu_1+\cdots +\nu_{n-1})
         \circ \cdots \circ  a^x_{\nu_2}(\tau+\nu_1) \circ a^x_{\nu_1}(\tau)$ 
     and obtain in this way an algebraic map
     \[
      \begin{array}{lccc}
       a :\, &\f{h}^* 
         & \longrightarrow 
         & \Hom_k \big(\f{U}(\f{n}^-), 
                  E(\nu_n)\otimes \cdots E(\nu_1) \otimes \f{U}(\f{n}^-)\big)\\
         & \tau   
         & \mapsto         
         & a^x_{\nu_n}(\tau+\nu_1+\cdots +\nu_{n-1}) \circ \cdots \circ
           a^x_{\nu_2}(\tau+\nu_1) \circ a^x_{\nu_1}(\tau)
       \end{array}
       \]      
 
      which vanishes only on a set of codimension $\ge 2$ and which maps a
      generic weight $\tau$ to
       $a(\tau) 
       = \delta_{x\nu_1}(x\cdot\tau)\delta_{x\nu_2}(x\cdot(\tau+\nu_1))
          \cdots \delta_{x\nu_n}(x\cdot(\tau+\nu_1+\cdots\nu_{n-1}))$
          $h^x(\nu_1, \ldots, \nu_n)(\tau)
       = d^x(\nu_1,\ldots,\nu_n)(\tau) h^x(\nu_1, \ldots, \nu_n)(\tau)$.
    We thus obtain the map $a$ as the desired algebraic extension.
  \end{proof}

  According to the Theorem of Bernstein-Gelfand we may interpret for generic
    $\tau$ the maps
    $h^x(\nu_1, \ldots, \nu_n)(\tau) \in 
    \Hom_{\f{g}} \big(\verma{x \cdot \tau},
               T(\nu_n)\cdots T(\nu_1)\verma{x \cdot \tau}\big)$ 
    $\cong \Homtau \big(\Id, T(\nu_n)\cdots T(\nu_1)\big)$ 
    as natural transformations of functors.
  Set now
   \begin{align*}
      h_1&:=h^x(\nu+\mu,-\mu,-\nu)(\tau)
       &&\in \Homtau\big(\Id, T(-\nu) T(-\mu) T(\nu+\mu)\big),\\
      h_2&:=h^x(-\nu,\nu)(\tau+\nu)
       &&\in \Hom_{{\cal M}(\chi(\tau+\nu))\to}\big(\Id,T(\nu) T(-\nu)\big),\\
      h_3&:=h^x(-\mu,\mu)(\tau+\nu+\mu)
       &&\in\Hom_{{\cal M}(\chi(\tau+\nu+\mu))\to}\big(\Id,T(\mu) T(-\mu)\big)
   \end{align*}
  
   and consider the natural transformations
   \[
   \begin{CD}
      T(\mu) T(\nu)
        @>{G^x(\tau)}>> 
      T(\nu+\mu)\\
      @V{T(\mu)T(\nu) h_1}VV
      @AA{(h_3)^{-1}}A \\
      T(\mu) T(\nu) T(-\nu) T(-\mu) T(\nu+\mu)
      @>{(T(\mu) h_2)^{-1}}>>
      T(\mu) T(-\mu) T(\nu+\mu)
   \end{CD}
   \]

   The diagram commutes, since  
    $(h_3)^{-1} \circ (T(\mu) h_2)^{-1}\circ T(\mu)T(\nu) h_1$ as well as 
    $G^x(\tau)$ imply the identity on $\verma{x\cdot(\tau+\nu+\mu)}$ under the
    isomorphism\\
\pagebreak

   \hspace*{0.5em}
    $\Homtau \big(T(\mu) \circ T(\nu), T(\nu+\mu)\big)$ \\[-2.5em]

    \begin{center}
    \begin{align*}
    &\stackrel{\sim}{\longrightarrow}
     \Hom_{\f{g}}\big(T(\mu) T(\nu)\verma{x\cdot \tau}, 
                                          T(\nu+\mu)\verma{x\cdot \tau}\big)\\ 
    &\stackrel{\sim}{\longrightarrow}
     \Hom_{\f{g}}\big(\verma{x\cdot(\tau+\nu+\mu)},
                                     \verma{x\cdot(\tau+\nu+\mu)}\big).
    \end{align*}
    \end{center}

   \begin{lemma}\label{delta}
      Set $D^x(\tau) := \delta_{x(\nu+\mu)}(x\cdot\tau)/
        \delta_{x\nu}(x\cdot\tau) \delta_{x\mu}(x\cdot(\tau+\nu))$.
      Then the map $D^x G^x$ is algebraic on $\cal U$ and there exists an
        algebraic extension on $\f{h}^*$ whose set of zeros has codimension
        $\ge 2$.
   \end{lemma}

   \begin{proof}
     Since the maps $d^x(\nu_1,\ldots,\nu_n) h^x(\nu_1,\ldots,\nu_n)$
        have such algebraic extensions (Lemma \ref{faktor}), the
        commutativity of the above diagram implies that also $D^x G^x$
        has such an algebraic extension, where
        \begin{align*}
        D^x(\tau) &= d^x(\nu+\mu,-\mu,-\nu)(\tau)d^x(-\nu,\nu)(\tau+\nu)
                                      d^x(-\mu,\mu)(\tau+\nu+\mu) \\
                  &= \delta_{x(\nu+\mu)}(x\cdot\tau)/
             \delta_{x\nu}(x\cdot\tau) \delta_{x\mu}(x\cdot(\tau+\nu)).
        \end{align*}
   \end{proof}
 
   Finally, we come to the 
   \begin{proof}[Proof of Theorem \ref{dreiecksfunktion}]
     By Lemma \ref{langeformel} it suffices to show that
        $\Delta(\mu,\nu;x)(\tau) = c \;\pi(\tau+\rho)$ for a non-vanishing
        constant $c$.
     Note, that $(D^x/D^e)(\tau) = \pm \pi(\tau+\rho)$.
     We deduce from Lemma \ref{gleichung} that for generic weights
        $D^x D^e G^e = \Delta(\mu,\nu;x) D^e D^x G^x$.
     Since $D^e G^e$ as well as $D^x G^x$ have algebraic extensions on 
        $\f{h}^*$ which vanish only on a set of codimension $\ge 2$ 
        (Lemma \ref{delta}), it follows that there is a constant 
        $c \in k^{\times}$, independent of $\tau, \mu, \nu$ and $x$, such that
        \[c \; \Delta(\mu,\nu;x)(\tau)=(D^x/D^e)(\tau) = \pm \pi(\tau+\rho).\]
   \end{proof}

     \section{\bf Outlook}

 \subsection{Identities}\label{zerl}
   There are many nice identities for the triangle functions.
   Obvious are the 

      \begin{norm} 
          \begin{align*}
              \Delta(\mu,\nu;e) &=1  &
              \Delta(0,\nu;x)   &=1  &
              \Delta(\nu,0;x)   &=1
          \end{align*}
      \end{norm}

    By means of Theorem \ref{dreiecksfunktion} or directly with the 
      definition of $\Delta$ one then checks for $\nu, \mu, \eta \in P$
      and $x, y \in \cal W$ : 

      \begin{zer}
               $\Delta(\eta+\mu,\nu;x)(\tau)\;\Delta(\eta,\mu;x)(\tau +\nu) 
              = \Delta(\eta,\mu+\nu;x)(\tau)\;\Delta(\mu,\nu;x)(\tau)$
      \end{zer}  

      \begin{rot}
           $\Delta(y\mu,y\nu;x)(y \cdot \tau)
            = (\Delta(\mu,\nu;y)(\tau))^{-1}\;\Delta(\mu,\nu;xy)(\tau)$
      \end{rot}

      \begin{flach}
           Let $\nu, \mu \in P$ be in the closure of a Weyl chamber. Then
           $\Delta(\mu,\nu;x) = 1$ for all $x \in \cal W$.
      \end{flach}

\pagebreak

 \subsection{Generalizations}
 \subsubsection{The Weyl group parameter}
    Let $\nu, \mu \in P$ and $x, y \in \cal W$.
    Instead of applying the translation functors to the Verma modules
       $\verma{\tau}$ and $\verma{x\cdot\tau}$, we may choose the Verma modules
       $\verma{x\cdot\tau}$ and $\verma{y\cdot\tau}$.
    We then define a generalized triangle function $\Deltag$ by
    \[\Deltag(\mu,\nu;y,x)(\tau) 
        := \det\left(yx^{-1} \circ \nat(\mu,\nu;x)(\tau) 
                               \circ (\nat(\mu,\nu;y)(\tau))^{-1} \right) .\]
    We now have
    $\Deltag(\mu,\nu;e,x)(\tau) = \Delta(\mu,\nu;x)(\tau)$
    and going back to the definition of $\Deltag$ we deduce the identities
      \begin{itemize}
         \item[{\bf (I1)}]$\Deltag(\mu,\nu;y,x)(\tau) 
                        = (\Delta(\mu,\nu;y)(\tau))^{-1} \;
                                              \Delta(\mu,\nu;x)(\tau)$

         \item[{\bf (I2)}]$\Deltag(\mu,\nu;y,x)(\tau) =
                         (\Deltag(\mu,\nu;x,y)(\tau))^{-1}$ 
      \end{itemize}
    and
      \begin{itemize} 
         \item[{\bf (I3)}] $\Deltag(\mu,\nu;y,x)\;\Deltag(\mu,\nu;x,z) 
                            = \Deltag(\mu,\nu;y,z)$.
      \end{itemize}
    
    The equivalent statement to the Rotation identity is obtained by
      comparing (I1) with

    \begin{itemize}
       \item[{\bf (I4)}] $\Deltag(\mu,\nu;y,x)(\tau) 
                   = \Delta(y\mu,y\nu;xy^{-1})(y\cdot\tau)$. 
    \end{itemize}

\subsubsection{Number of translations}

   The triangle functions measure in a subtle way the relation between the
     two translation functors 
     $T_{\tau+\nu}^{\tau+\nu+\mu} \circ T_{\tau}^{\tau+\nu}$ and 
     $T_{\tau}^{\tau+\nu+\mu}$. 
   Therefore the triangle

       {\xymatrix{
        &&&& & {\centerdot} \ar[dr]^{\mu}\\
        &&&&{\tau \centerdot} \ar[ur]^{\nu} \ar[rr]_{\nu+\mu} && {\centerdot}}}

   Let now integral weights $\nu_{1}, \ldots ,\nu_{n}$ and 
     $\mu_{1}, \ldots ,\mu_{m} \in P$ be given such that
     $\sum_{i=1}^n \nu_i = \sum_{j=1}^m \mu_j =: p$.
   Call then $T(\nu_i)$ the translation functor
     {\renewcommand{\arraystretch}{1.5}
      \[
      \begin{array}{lccc}
       T(\nu_i):\, 
         & \cal M^{\infty} \big(\chi(\tau+\nu_1+\cdots+\nu_{i-1})\big) 
         & \longrightarrow 
         & \cal M^{\infty} \big(\chi(\tau+\nu_1+\cdots+\nu_i)\big)\\
         & M   
         & \mapsto 
         & \pr_{\chi(\tau+\nu_1+\cdots+\nu_i)}(E(\nu_i) \otimes M)
     \end{array}
     \]}

     and define similarly the translation functor $T(\mu_i)$.
 We now want to compare the functors 
     $T(\nu_n)\circ \cdots \circ T(\nu_1)$ 
     and 
     $T(\mu_m)\circ \cdots \circ T(\mu_1)$ 
     with each other and instead of a triangle of translations we thus
     have now the following situation:

     {\xymatrix{
        &&&&{\centerdot} \ar[r]
        &{\centerdot} \ar@{.}[r]
        &{\centerdot} \ar[rd]^{\mu_m}\\
        &&&{\tau\centerdot} \ar[ur]^{\mu_1} \ar[dr]_{\nu_1}&&&&
                                               {\centerdot \tau+p}\\
        &&&&{\centerdot} \ar[r]
        &{\centerdot} \ar@{.}[r]
        &{\centerdot} \ar[ru]_{\nu_n}}}

  We start by defining a map 
     $\bar{\nat}(\mu_m,\ldots,\mu_1;\nu_n,\ldots,\nu_1;x)(\tau)$ for
     $x\in \cal W$ and generic weight $\tau$ analogously to the definition
     of $\nat(\mu,\nu;x)(\tau)$ (see page \pageref{nat}) as the composition \\

\pagebreak

   \hspace*{0.5em}
   $\Homtau \big(T(\mu_m)\circ \cdots \circ T(\mu_1),
            T(\nu_n)\circ \cdots \circ T(\nu_1)\big)$\\[-2.7em]

    \begin{center}
    \begin{align*}
    &\stackrel{\sim}{\longrightarrow} 
     \Hom_{\f{g}}\big(T(\mu_m) \cdots T(\mu_1)\verma{x \cdot \tau},
                  T(\nu_n) \cdots T(\nu_1)\verma{x \cdot \tau}\big)\\
    &\stackrel{\sim}{\longrightarrow} 
     \Hom_{\f{g}}\big(E_{\mu_m} \hotimes \cdots E_{\mu_1} \hotimes 
                                           \verma{x \cdot (\tau +p)},
       E_{\nu_n}\hotimes\cdots E_{\nu_1}\hotimes\verma{x\cdot(\tau+p)}\big)\\
    &\stackrel{\sim}{\longrightarrow} 
           E_{\mu_m}^* \otimes \cdots \otimes E_{\mu_1}^* 
                     \otimes E_{\nu_n} \otimes \cdots \otimes E_{\nu_1}
    \end{align*}
    \end{center}

   Here, we wrote $E_{\mu}$ for $E(\mu)_{x\mu}$.
   We then obtain a generalized triangle function $\bDelta$ by
    \begin{multline*}
    \bDelta(\mu_m,\ldots\mu_1;\nu_n,\ldots\nu_1;x)(\tau)
       := \det \big(x^{-1}\circ 
       \bar{\nat}(\mu_m,\ldots\mu_1;\nu_n,\ldots\nu_1;x)(\tau)\\ \circ 
       (\bar{\nat}(\mu_m,\ldots\mu_1;\nu_n,\ldots\nu_1;e)(\tau))^{-1}\big).
     \end{multline*} 

   Obviously we have
     $ \bDelta(\mu,\nu;\nu+\mu;x) = \Delta(\mu,\nu;x)$ as well as
     \[\bDelta(\mu_m,\ldots,\mu_1;\nu_n,\ldots,\nu_1;x)(\tau) 
           = (\bDelta(\nu_n,\ldots,\nu_1;\mu_m,\ldots,\mu_1;x)(\tau))^{-1}\]

     and we can reduce the calculation of $\bDelta$ to the calculation of 
     $\Delta$ by means of the

   \begin{Split}
       \begin{multline}
       \bDelta(\nu_3,\nu_2,\nu_1;\nu_3+\nu_2+\nu_1;x)(\tau)
       = \Delta(\nu_2,\nu_1;x)(\tau)\Delta(\nu_3,\nu_1+\nu_2;x)(\tau)\\
       = \Delta(\nu_2+\nu_3,\nu_1;x)(\tau)\Delta(\nu_3,\nu_2;x)(\tau+\nu_1)
      \end{multline}
       \begin{equation}
       \bDelta(\mu_2,\mu_1;\nu_2,\nu_1;x)(\tau)
       =\Delta(\mu_2,\mu_1;x)(\tau)\big(\Delta(\nu_2,\nu_1;x)(\tau)\big)^{-1}
      \end{equation}
   \end{Split}

  The first identity just means that we can split \hspace{0.5em}
      \begin{minipage}{4.5cm}
         {\xymatrix{
         {\centerdot}\ar[rr]^{\nu_2} 
         &&{\centerdot} \ar[d]^{\nu_3}\\
         {\tau\centerdot} \ar[u]^{\nu_1} \ar[rr]_{\nu_1+\nu_2+\nu_3}
         &&{\centerdot}}}
         \end{minipage}\hspace{0.5em}
    into triangles according to \hspace{0.5em}
       \begin{minipage}{4.5cm}
         {\xymatrix{
         {\centerdot} \ar[rr]^{\nu_2}
         &&{\centerdot} \ar[d]^{\nu_3} \\
         {\tau\centerdot} \ar[u]^{\nu_1} \ar@<-0.5ex>[rr]_{\nu_1+\nu_2+\nu_3}
                     \ar@<0.3ex>[urr]\ar@<-0.7ex>[urr]
         &&{\centerdot}}}
       \end{minipage}\hspace{0.5em}
   or \hspace{0.5em}
       \begin{minipage}{4.5cm}
         {\xymatrix{
         {\centerdot} \ar[rr]^{\nu_2}\ar@<0.3ex>[rrd]\ar@<-0.7ex>[rrd]
         &&{\centerdot} \ar[d]^{\nu_3} \\
         {\tau\centerdot} \ar[u]^{\nu_1} \ar@<-0.5ex>[rr]_{\nu_1+\nu_2+\nu_3}
         &&{\centerdot}}}
       \end{minipage}.
 This is just the Decomposition identity in \ref{zerl}.
 The second equation describes the decomposition of \hspace{0.5em}
     \begin{minipage}{5.5cm}
         {\xymatrix{
            &{\centerdot}\ar[dr]^{\mu_2}\\
            {\tau\centerdot}\ar[ur]^{\mu_1}\ar[dr]_{\nu_1}&&{\centerdot}\\
            &{\centerdot}\ar[ur]_{\nu_2}}}
      \end{minipage}\hspace{0.5em}
   in \hspace{0.5em} 
      \begin{minipage}{5.5cm}
        {\xymatrix{
          &{\centerdot}\ar[dr]^{\mu_2}\\
          {\tau\centerdot}\ar[ur]^{\mu_1}\ar[dr]_{\nu_1}
                           \ar@<0.5ex>[rr] \ar@<-0.5ex>[rr]&&{\centerdot}\\
          &{\centerdot}\ar[ur]_{\nu_2}}} 
      \end{minipage}

  Inductively, we may thus first split up our situation into

       
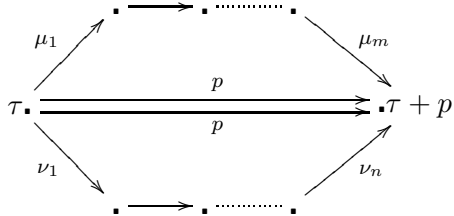
{\xymatrix{
        &&&&{\centerdot} \ar[r]
        &{\centerdot} \ar@{.}[r]
        &{\centerdot} \ar[rd]^{\mu_m}\\
        &&&{\tau \centerdot} \ar[ur]^{\mu_1} \ar[dr]_{\nu_1}
                \ar@<0.5ex>[rrrr]^{p}
                \ar@<-0.5ex>[rrrr]_{p}
                        &&&&{\centerdot\tau+p}\\
        &&&&{\centerdot} \ar[r]
        &{\centerdot} \ar@{.}[r]
        &{\centerdot} \ar[ru]_{\nu_n}}}

  In formulae :
     \begin{align*}
     \bDelta(\mu_m,\ldots \mu_1;\nu_n,\ldots \nu_1;x)(\tau)
          &=\bDelta(\mu_m,\ldots \mu_1;p;x)(\tau)
           \bDelta(p;\nu_n,\ldots \nu_1;x)(\tau)\\
          &=\bDelta(\mu_m,\ldots \mu_1;p;x)(\tau)
           \big(\bDelta(\nu_n,\ldots \nu_1;p;x)(\tau)\big)^{-1}
    \end{align*}

     and in order to calculate the $\bDelta$-functions it suffices thus to
     know them in the special case
     $\bDelta(\mu_m,\ldots , \mu_1;\sum_{j=1}^m \mu_j;x)$.
  This situation can then be reduced to triangles by decomposing it into

        
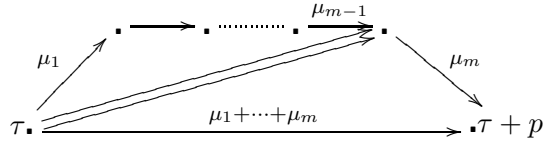
{\xymatrix{
        &&&{\centerdot} \ar[r]
        &{\centerdot} \ar@{.}[r]
        &{\centerdot} \ar[r]^{\mu_{m-1}}
        &{\centerdot} \ar[rd]^{\mu_m}\\ 
        &&{\tau \centerdot} \ar[ur]^{\mu_1}
                          \ar[urrrr]
                          \ar@<-0.7ex>[urrrr]
                          \ar@<-0.5ex>[rrrrr]^{\mu_1+\cdots+\mu_m}
                        &&&&&{\centerdot\tau+p}}}

   or into

        
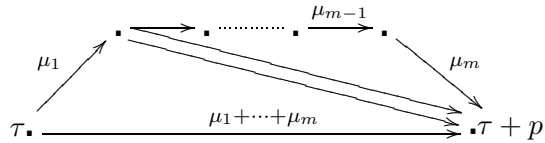
{\xymatrix{
        &&&{\centerdot} \ar[r]\ar@<0.3ex>[drrrr]\ar@<-0.7ex>[drrrr]
        &{\centerdot} \ar@{.}[r]
        &{\centerdot} \ar[r]^{\mu_{m-1}}
        &{\centerdot} \ar[rd]^{\mu_m}\\ 
        &&{\tau \centerdot} \ar[ur]^{\mu_1}
                          \ar@<-0.5ex>[rrrrr]^{\mu_1+\cdots+\mu_m}
                        &&&&&{\centerdot\tau+p}}} 
 
   Inductively one can then prove
  
   \begin{proposition}
     Let $\mu_1, \ldots , \mu_m \in P$ be integral weights and $x \in \cal W$.
     Then
     \begin{align*}
     \bDelta(\mu_m,\ldots,\mu_1;\textstyle \sum\nolimits_{j=1}^m \mu_j;x)(\tau)
        &=\prod_{k=2}^m 
             \Delta(\mu_k,\textstyle \sum\nolimits_{i=1}^{k-1}\mu_i;x)(\tau)\\
        &=\prod_{k=1}^{m-1} 
             \Delta(\textstyle \sum\nolimits_{j=k+1}^m \mu_j,\mu_k;x)
                 (\tau+\textstyle \sum\nolimits_{j=1}^{k-1}\mu_j)
        \end{align*}
   \end{proposition}



\end{document}